\documentclass[a4paper,11pt]{article}

\usepackage{graphicx,psfrag}
\usepackage{amssymb,amsmath,amsthm}
\usepackage{enumerate}  
\usepackage{a4wide}
\usepackage{blindtext}
\usepackage{caption}
\usepackage{subcaption}

\usepackage{mathrsfs}
\usepackage{enumitem}

\numberwithin{equation}{section}
\allowdisplaybreaks[4]

\newtheorem{proposition}{Proposition}[section]
\newtheorem{theorem}[proposition]{Theorem}
\newtheorem{lemma}[proposition]{Lemma}
\newtheorem{corollary}[proposition]{Corollary}
\newtheorem{definition}[proposition]{Definition}
\theoremstyle{definition}
\newtheorem{remark}[proposition]{Remark}
\newtheorem{example}{Example}
\newtheorem*{ack}{Acknowledgments}

\begin{document}

\title{Non-linear functionals, deficient topological measures, and representation theorems on locally compact spaces}
\author{Svetlana V. Butler\footnote{Department of Mathematics,  University of California, Santa Barbara,  552 University Rd, Isla Vista, CA 93117, USA. 
Email: svtbutler@ucsb.edu}
}
\maketitle

\maketitle
\begin{abstract}
We study non-linear functionals, including quasi-linear functionals, 
p-conic quasi-linear functionals, d-functionals, r-functionals, and their relationships to deficient topological measures and 
topological measures on locally compact spaces. We prove representation theorems and
show, in particular, that there is an order-preserving, conic-linear bijection between the class of finite deficient topological measures
and the class of bounded p-conic quasi-linear functionals. 
Our results imply known representation theorems for finite topological measures and  deficient topological measures.
When the space is compact we obtain four equivalent definitions of a quasi-linear functional and four 
equivalent definitions of functionals corresponding to deficient topological measures.
\end{abstract}



\section{Introduction}

This paper is devoted to the study of various non-linear functionals and their relationships to deficient topological measures and  
topological measures on locally compact spaces.  
J. F.  Aarnes first discovered quasi-linear functionals and the corresponding set functions, topological measures, 
(initially called quasi-measures) in~\cite{Aarnes:TheFirstPaper}. 
Since then many works devoted to these objects have appeared. For their use in symplectic topology one may consult 
numerous papers beginning with~\cite{EntovPolterovich} and a monograph (\cite{PoltRosenBook}). 

We prove representation theorems for deficient topological measures on locally compact spaces. 
Our results imply known results, including the representation theorem for quasi-linear functionals 
on compact spaces (\cite{Aarnes:TheFirstPaper}),  
the representation theorem for deficient topological measures on compact spaces (\cite{Svistula:DTM}), and 
the representation theorems for quasi-linear functionals on locally compact spaces (\cite{Alf:ReprTh},\cite{Butler:QLFLC}). 
We also obtain new consequences where $X$ is compact, including four equivalent definitions of a quasi-linear functional and four 
equivalent definitions of functionals corresponding to deficient topological measures.
 
Deficient topological measures were first defined and used by A. Rustad 
and $\emptyset$. Johansen in~\cite{OrjanAlf:CostrPropQlf}. 
They were later independently rediscovered  by M. Svistula (\cite{Svistula:Signed},~\cite{Svistula:DTM}).  
In these works, deficient topological measures were defined as real-valued functions on a compact space. In this paper
we use deficient topological measures on locally compact spaces as functions into extended real numbers. Quasi-linear functionals
have been generalized to other non-linear functionals in order to represent deficient topological measures on compact spaces.
See~\cite[section 6]{OrjanAlf:CostrPropQlf};  see also~\cite{Svistula:DTM}, where such functionals are called r- and l-functionals.
To represent deficient topological measures on locally compact spaces we will need the above mentioned non-linear functionals 
and some others. This paper is influenced by many works, 
including~\cite{Aarnes:TheFirstPaper},~\cite{OrjanAlf:CostrPropQlf},~\cite{Svistula:DTM},~\cite{Alf:ReprTh},~\cite{Grubb:Signed}, and~\cite{Butler:QLFLC}.
 
The paper is organized as follows. Section \ref{SePrelim} contains necessary definitions and background facts. 
In Section \ref{SeCones} we consider cones generated by a function, as well as p-conic and n-conic quasi-linear functionals.
In Section \ref{SeDfnl} we consider various related non-linear functionals, including d-, r-, l-, and s-functionals. 
In Section \ref{SeGetDTM}  we show how to obtain deficient topological measures from d-functionals.
In Section \ref{lerim}, given a finite deficient topological measure and a bounded continuous function, 
we define left and right measures. We give a criterion for left and right measures to coincide. Left and right measures coincide if 
one starts from a finite topological measure; 
we also give examples which show that left and right measures 
may or may not coincide even for $\{0,1\}$-valued deficient topological measures. 
In section \ref{SeGetFnls} using right and left measures we obtain by integration p-conic and l-conic quasi-linear functionals, 
which are also r- and l-functionals.
Section \ref{SeReprT} is devoted to representation theorems for finite deficient topological measures in terms of bounded 
p-conic and  l-conic quasi-linear functionals, and also r-  and l-functionals. In particular, we show that 
there is an order-preserving, conic-linear bijection between the class of finite deficient topological measures
and  the class of bounded p-conic quasi-linear functionals. 
When $X$ is compact, we obtain four equivalent definitions of quasi-linear functionals and four 
equivalent definitions of functionals corresponding to deficient topological measures.

\section{Preliminaries} \label{SePrelim}

In this paper $X$ is a locally compact, connected space. 

By $C(X)$ we denote the set of all real-valued continuous functions on $X$ with the uniform norm, 
by $C_b(X)$ the set of bounded continuous functions on $X$,   
by $C_0(X)$ the set of continuous functions on $X$ vanishing at infinity,  and
by $C_c(X)$ the set of continuous functions with compact support. 
By $C_0^+(X)$ we denote the collection of all nonnegative functions vanishing at infinity; 
similarly, $C_0^-(X)$ is the collection of all nonpositive
functions from $C_0(X)$.

When we consider maps  into $ [-\infty, \infty]$ we assume that any such map 
attains at most one of $ \infty, - \infty$, and is not identically $\infty$ or $ - \infty$. 
By $D(\rho)$ we denote the domain of a functional $\rho$. For example, we may take $D(\rho) = C_0^+(X)$ 
or  $C_c(X)$. If $ D(\rho) = C(X)$, where $X$ is compact, then $D(\rho)$ contains constants.

We denote by $\overline E$ the closure of a set $E$, and by $ \bigsqcup$ a union of disjoint sets.
We denote by $1$ the constant function $1(x) =1$,  by $id$ the identity function $id(x) = x$, 
and by $1_K$ the characteristic function of a set $K$. By $ supp \,  f $ we mean $ \overline{ \{x: f(x) \neq 0 \} }$.

Several collections of sets are used often.   They include:
$\mathscr{O}(X)$,  the collection of open subsets of   $X $;
$\mathscr{C}(X)$  the collection of closed subsets of   $X $;
$\mathscr{K}(X)$  the collection of compact subsets of   $X $.

\begin{definition} \label{MDe2}
Let $X$ be a  topological space and $\mu$ be a nonnegative set function on $\mathcal{E}$, 
a family of subsets of $X$ that contains $\mathscr{O}(X) \cup \mathscr{C}(X)$. 
We say that 
\begin{itemize}
\item
$\mu$ is inner regular (or inner compact regular) 
if $\mu(U) = \sup \{ \mu(C) : \ \ C \subseteq U, \ \ C \in  \mathscr{K}(X) \}$ for each open set $U$.
\item
$\mu$ is outer regular if 
$\mu(F) = \inf \{ \mu(U) : \ \ F \subseteq U, \ \ U \in  \mathscr{O}(X) \}$ for each closed set $F$. 
\item
$ \parallel \mu \parallel = \mu(X)$ and $\mu$ is finite if $\mu(X) < \infty$.
\item
$\mu$ is compact-finite if $\mu(K) < \infty$ for any $ K \in \mathscr{K}(X)$.
\item
$ \mu$ is monotone if $ A \subseteq B$ implies $ \mu(A) \le \mu(B)$.
\item
$\mu$ is $\tau$-smooth on compact sets if for every decreasing net
$K_\alpha \searrow K, K_\alpha, K \in \mathscr{K}(X)$ we have $\mu(K_\alpha) \rightarrow \mu(K)$.
\item
$\mu$ is $\tau$-smooth on open sets if  for every increasing net  
$U_\alpha \nearrow U, U_\alpha, U \in \mathscr{O}(X)$ we have $\mu(U_\alpha) \rightarrow \mu(U)$.
\item
$\mu$ is simple if it only assumes  values $0$ and $1$.
\end{itemize} 
\end{definition}

\begin{definition}\label{DTM}
A  deficient topological measure on a locally compact space $X$ is a set function
$\nu:  \mathscr{C}(X) \cup \mathscr{O}(X) \longrightarrow [0, \infty]$ 
which is finitely additive on compact sets, inner compact regular, and 
outer regular, i.e. :
\begin{enumerate}[label=(DTM\arabic*),ref=(DTM\arabic*)]
\item \label{DTM1}
if $C \cap K = \O, \ C,K \in \mathscr{K}(X)$ then $\nu(C \sqcup K) = \nu(C) + \nu(K)$; 
\item \label {DTM2} 
$ \nu(U) = \sup\{ \nu(C) : \ C \subseteq U, \ C \in \mathscr{K}(X) \} $
 for $U\in\mathscr{O}(X)$;
\item \label{DTM3} 
$ \nu(F) = \inf\{ \nu(U) : \ F \subseteq U, \ U \in \mathscr{O}(X) \} $  for  $F \in \mathscr{C}(X)$.
\end{enumerate}
\end{definition}
\noindent
For a closed set $F$, $ \nu(F) = \infty$ iff $ \nu(U) = \infty$ for every open set $U$ containing $F$.

\begin{remark} \label{DTMagree}
Note that a deficient topological measure $ \nu$ is monotone on $ \mathscr{O}(X) \cup \mathscr{C}(X)$ and $\nu(\O) = 0$.
If $\nu$ and $\mu$ are deficient topological measures that agree on $\mathscr{K}(X)$ then $\nu =\mu$, and 
if $\nu \le \mu$ on $\mathscr{K}(X)$ (or on $ \mathscr{O}(X)$) then $\nu  \le \mu$.
\end{remark}

\begin{definition}\label{TMLC}
A topological measure on $X$ is a set function
$\mu:  \mathscr{C}(X) \cup \mathscr{O}(X) \to [0,\infty]$ satisfying the following conditions:
\begin{enumerate}[label=(TM\arabic*),ref=(TM\arabic*)]
\item \label{TM1} 
if $A,B, A \sqcup B \in \mathscr{K}(X) \cup \mathscr{O}(X) $ then
$
\mu(A\sqcup B)=\mu(A)+\mu(B);
$
\item \label{TM2}  
$
\mu(U)=\sup\{\mu(K):K \in \mathscr{K}(X), \  K \subseteq U\}
$ for $U\in\mathscr{O}(X)$;
\item \label{TM3}
$
\mu(F)=\inf\{\mu(U):U \in \mathscr{O}(X), \ F \subseteq U\}
$ for  $F \in \mathscr{C}(X)$.
\end{enumerate}
\end{definition} 

We denote by $TM(X)$ the collection of all topological measures on $X$, 
by $DTM(X)$ the collection of all deficient topological measures on $X$, and by
$M(X)$ the collection of all Borel measures on $X$ that are inner regular on open sets and outer regular
(restricted to $\mathscr{O}(X) \cup \mathscr{C}(X)$). 

\begin{remark} \label{proinclu}
Let $X$ be locally compact. In general,
$$ M(X) \subsetneqq  TM(X) \subsetneqq  DTM(X). $$
When $X$ is compact, there are examples of topological measures that are not measures 
and of deficient topological measures that are not topological measures in numerous papers, 
beginning with~\cite{Aarnes:TheFirstPaper},~\cite{OrjanAlf:CostrPropQlf}, and~\cite {Svistula:Signed}.
When $X$ is locally compact, see~\cite{Butler:TechniqLC},~\cite[Sections 5 and 6]{Butler:DTMLC}, and~\cite[Section 9]{Butler:TMLCconstr} 
for more information on proper inclusion, criteria for a deficient topological measure to belong to $M(X)$ or $TM(X)$ (in particular, to be a Radon measure 
or a regular Borel measure), as well as various examples.
\end{remark}

The proof of the next result is in~\cite[Section 4]{Butler:DTMLC}.
\begin{theorem} \label{DTMtoTM}
\begin{enumerate}[label=(\Roman*),ref=(\Roman*)]
\item
Let $X$ be compact, and $\nu$ a deficient topological measure. The following are equivalent:
\begin{enumerate}
\item[(a)]
$\nu$ is a topological measure.
\item[(b)]
$\nu(X) = \nu(C)  + \nu(X \setminus C), \ \ \ C \in \mathscr{C}(X).$ 
\item[(c)]
$\nu(X) \le \nu(C)  + \nu(X \setminus C), \ \ \ C \in \mathscr{C}(X).$  
\end{enumerate}
\item
Let $X$ be locally compact, and $\nu$ a deficient topological measure. 
The following are equivalent:
\begin{enumerate}
\item[(a)]
$\nu$ is a topological measure.
\item[(b)]
$\nu(U) = \nu(C)  + \nu(U \setminus C), \ \ \ C \in \mathscr{K}(X), \ \ U \in \mathscr{O}(X).$ 
\item[(c)]
$\nu(U) \le \nu(C)  + \nu(U \setminus C), \ \ \ C \in \mathscr{K}(X), \ \ U \in \mathscr{O}(X).$
\end{enumerate}
\end{enumerate}
\end{theorem}

Recall the following fact (see, for example,~\cite[Chapter XI, 6.2]{Dugundji}):
\begin{lemma} \label{easyLeLC}
Let $K \subseteq U, \ K \in \mathscr{K}(X),  \ U \in \mathscr{O}(X)$ in a locally compact space $X$.
Then there exists a set  $V \in \mathscr{O}(X)$ with compact closure such that
$$ K \subseteq V \subseteq \overline V \subseteq U. $$ 
\end{lemma}

The following is proved in~\cite[Section 3]{Butler:DTMLC}.
\begin{lemma} \label{opaddDTM}
Let $X$ be a locally compact space.
\begin{itemize}
\item[(a)]
A deficient topological measure is $\tau$-smooth on compact sets and 
$\tau$-smooth on open sets. In particular, a topological measures is additive on open sets. 
\item[(b)]
A deficient topological measure $ \mu$ is superadditive, i.e. 
if $ \bigsqcup_{t \in T} A_t \subseteq A, $  where $A_t, A \in \mathscr{O}(X) \cup \mathscr{C}(X)$,  
and at most one of the closed sets (if there are any) is not compact, then 
$\mu(A) \ge \sum_{t \in T } \mu(A_t)$. 
\end{itemize}
\end{lemma}

We recall the following theorem:

\begin{theorem}[Integration by parts for Lebesque-Stieltjes integrals] \label{LebSt}
Let $ F$ and $ G$ be real-valued nondecreasing functions on $\mathbb{R}$ with corresponding Lebesque-Stieltjes measures
$m_{F}, \ m_{G}$. Then for $ a < b$ we have:
$$ \int_a^b G(x^+) dm_{F}(x) + \int_a^b F(x^-) dm_{G}(x) = F(b^+) G(b^+) - F(a^-) G(a^-).$$
\end{theorem}

\begin{definition} \label{fnls}
Let $\rho$ be a functional on $C(X)$. We say   
\begin{itemize}
\item
$\rho$ is homogeneous if $\rho (a f) = a \rho(f)$  for any $ a \in \mathbb{R}$.
\item
$\rho$ is positive-homogeneous if  $\rho (a f) = a \rho(f)$  for any $ a \ge 0$.
\item
$\rho$ is conic-linear on a cone $S$ if it is linear on conic combinations of elements of $S,$ 
i.e. $\rho (a f + b g) = a \rho(f) + b \rho(g)$  for any $ a,b \ge 0$ and any $f , g \in S$. 
\item
$\rho$ is positive if $f \ge 0 \Longrightarrow  \rho(f) \ge 0$.
\item
$\rho$ is monotone if $ g \le f \Longrightarrow  \rho(g) \le \rho(f) $.
\item
$\rho$ is orthogonally additive if  $f \cdot g = 0  \Longrightarrow  \rho(f + g) = \rho(f) + \rho(g)$.
\item
$\rho$ is real-valued if $ \rho(f)  \in \mathbb{R}$ for any $f$. 
\item
$\parallel \rho \parallel =  \sup \{ | \rho(f) | : \  \parallel f \parallel \le 1 \} $ and $\rho$ is bounded if $\parallel \rho \parallel < \infty$.
\end{itemize}
\end{definition}

The  remaining definitions and facts related to quasi-linear functionals can be found in~\cite[Section 2]{Butler:QLFLC}:

\begin{definition} \label{sigesuba}
Let $X$  be locally compact.
\begin{itemize}
\item[(a)] 
Let $f \in C_b(X)$. Define $A(f)$  to be the smallest closed subalgebra 
of $C_b(X)$ containing $f$ and $1$. 
Hence, when $X$ is compact, we take  $f \in C(X)$ and define $A(f)$  
to be the smallest closed subalgebra of $C(X)$ containing $f$ and $1$. 
We call $A(f)$ the singly generated subalgebra  of $C(X)$ generated by $f$.
\item[(b)]
Let  $\mathcal{B}$ be a sublagebra of  $C_b(X)$. Define $B(f)$ to be the smallest closed subalgebra 
of $\mathcal{B}$ containing $f$. We call $B(f)$ the singly generated subalgebra  
of $\mathcal{B}$ generated by $f$. 
\end{itemize}
\end{definition}

We may take, for example,  $C_c(X), \ C_0(X)$ as $\mathcal{B}$.   

\begin{remark} \label{smsubalg}
When $X$ is compact,  $A(f)$ for $f \in C(X)$ contains all polynomials of $f$. It is not hard to show that $A(f)$ has the form:
\[ A(f)  = \{ \phi \circ f :  \phi \in C({f(X)}) \}. \]
When $X$ is locally compact, $\mathcal{B} = C_0(X)$ and $f \in C_0(X)$  (or $\mathcal{B} = C_c(X)$ and  $f \in C_c(X)$)  
the singly generated subalgebra  has the form:
\[ B(f) =  \{ \phi \circ f :  \phi(0) = 0, \ \phi \in C(\overline{f(X)}) \}. \]  
\end{remark}

\begin{definition} \label{QI}
Let $X$ be locally compact, and let $\mathcal{B}$ be a subalgebra of $C(X)$
containing $C_c(X)$.
A real-valued map $\rho$ on $ \mathcal{B}$ 
is a signed quasi-linear functional on $\mathcal{B}$
if 
\begin{enumerate}[label=(QI\arabic*),ref=(QI\arabic*)]
\item \label{QIconsLC}
$\rho(a f) = a \rho(f)$ for $ a \in \mathbb{R};$
\item \label{QIlinLC}
for each  $h \in \mathcal{B}  $ we have:
$\rho(f + g) =  \rho (f) + \rho (g)$ for $f,g$ in the singly generated subalgebra $B(h)$ generated by $h$. \\

We say that $\rho$ is a quasi-linear functional (or a positive quasi-linear functional)  if, in addition, 
\item \label{QIpositLC}
$ f \ge 0 \Longrightarrow \rho(f)  \ge 0.$
\end{enumerate}
When $X$ is compact, we call $\rho$ a quasi-state if $\rho(1) =1$.
\end{definition}

\begin{definition}
We say that a signed quasi-linear functional $\rho$ is compact-finite
if $| \rho(f)| < \infty$ for $f \in C_c(X)$.
\end{definition}
 
\section{Cones generated by functions} \label{SeCones}

When $X$ is  a compact or locally compact, non-compact space,
there is a correspondence between topological measures and quasi-linear functionals. 
(See~\cite{Aarnes:TheFirstPaper},~\cite{Alf:ReprTh},~\cite{Butler:QLFLC}). 
When $X$ is compact, there is also a correspondence between deficient topological measures and r- or l-functionals,
see~\cite{Svistula:DTM} and~\cite {OrjanAlf:CostrPropQlf}.
When $X$ is locally compact, we shall consider relations 
between deficient topological measures and various non-linear functionals. We shall start with 
functionals that are conic-linear on cones generated by functions.  

\begin{definition} \label{A+}
Let $X$ be compact, $f \in C(X)$. 
Define cones: 
$$ A^+(f) = \{ \phi \circ f: \ \phi  \mbox{   is non-decreasing   } \} ,$$
$$ A^-(f) = \{ \phi \circ f: \ \phi \mbox{   is non-increasing   } \}. $$

If $X$ is  a locally compact, non-compact space, $f  \in C_0(X)$
let 
$$ A^+(f) = \{ \phi \circ f: \ \phi  \mbox{   is non-decreasing,   } \phi(0) = 0 \} ,$$
$$ A^-(f) = \{ \phi \circ f: \ \phi \mbox{   is non-increasing,   } \phi(0) = 0  \}. $$
\end{definition}

\begin{remark} \label{fA+}
Since $ f = id \circ f$, for $ f \neq 0 $ ($ f \neq const $ if $X$ is compact) we have $ f \in A^+(f), f \notin  A^-(f)$, and $ -f \in A^-(f), -f \notin   A^+(f)$.
Note that $0 \in   A^+(f),  A^-(f)$. Also, $f^+ \in  A^+(f), f^- \in  A^-(f)$, 
since $ f^+ =  (id \vee 0) \circ f,  f^- = ( (-id) \vee 0) \circ f$.
Obviously, $  A^+(f),  A^-(f) \subseteq B(f)$ (respectively,  $  A^+(f),  A^-(f) \subseteq A(f)$ if $X$ is compact).
\end{remark} 

\begin{definition} \label{cqlf}
We call a  functional $\rho$ on $C_0(X)$  with values in $[ -\infty, \infty]$ (assuming at most one of $\infty, - \infty$) 
and $| \rho(0) | < \infty$ a p-conic quasi-linear functional if 
it is orthogonally additive and monotone on nonnegative functions and 
conic-linear on $A^+(f)$ for each $ f \in  C_0(X)$, i.e. 
\begin{enumerate}[label=(p\arabic*),ref=(p\arabic*)]
\item
If $f\, g=0, f, g  \ge 0$ then $ \rho(f+ g) = \rho(f) + \rho(g)$.
\item
If $0 \le g \le f$ then $\rho(g) \le \rho(f)$.
\item
For each $f$, if $g,h \in A^+(f), \ a,b \ge 0$ then $\rho(a g + bh) = a \rho(g) + b \rho(h)$.
\end{enumerate}

Similarly, a functional $\rho$ on $C_0(X)$  with values in $[ -\infty, \infty]$ (assuming at most one of $\infty, - \infty$) 
and $| \rho(0) | < \infty$ is called an n-conic quasi-linear functional 
if it is orthogonally additive and monotone on non-positive functions and 
conic-linear on $A^-(f)$ for each $ f \in  C_0(X)$. In other words, 
\begin{enumerate}[label=(n\arabic*),ref=(n\arabic*)]
\item
If $f\, g=0, f, g  \le 0$ then $ \rho(f+ g) = \rho(f) + \rho(g)$.
\item 
If $f \le g \le 0 $ then $\rho(f) \le \rho(g)$.
\item
For each $f$, if $g,h \in A^-(f), \ a,b \ge 0$ then $\rho(a g + bh) = a \rho(g) + b \rho(h)$.
\end{enumerate}
\end{definition}  

\begin{remark} \label{QIpnconic}
From~\cite[Lemma 20(iii), Section 3 and Lemma 44, Section 4]{Butler:QLFLC} it follows that a positive quasi-linear functional is 
a p-conic quasi-linear functional and also an n-conic quasi-linear functional.
\end{remark} 

\begin{remark} \label{PiA+}
Given a functional $\rho$, consider also the functional $\pi$ defined by $\pi(f) = - \rho(-f)$ for every $f \in D(\rho)$.  
Then $\pi$ is an n-conic quasi-linear functional iff
$ \rho$ is a p-conic quasi-linear functional. This allows us to transfer results for 
p-conic quasi-linear functionals to n-conic quasi-linear functionals and vice versa.
\end{remark}

\begin{lemma} \label{cqiConst} 
Let $ \rho$ be a functional on a locally compact space with $| \rho(0) | < \infty$.
\begin{enumerate}[label=(\Roman*),ref=(\Roman*)]
\item
Suppose a functional $\rho$ is conic-linear on $A^+(f)$ or on $A^-(f) $ for some $f$.
Then $ \rho(0) = 0.$
\item  \label{addPO}
Suppose a functional $\rho$ is conic-linear on $A^+(h)$ for each function $h \in C_0(X)$. Then  $ fg=0, f \ge 0, g \le 0 $ 
implies  $\rho(f+g) = \rho(f)  + \rho(g)$.
\item \label{levcoL}
Suppose $0 \le g(x) \le f(x) \le c, f=c \mbox{   on  } \{ x: g(x) >0\} $. Then $f,g \in A^+(f+g)$. 
In particular, if $ \rho$ is a p-conic quasi-linear functional then $ \rho(af+bg) =a\rho(f) + b\rho(g)$ for any $a, b \ge 0$.
If $ \rho$ is a quasi-linear functional, then $ \rho(af+bg) =a\rho(f) + b\rho(g)$ for any $a, b \in \mathbb{R}$.
\item \label{levcoMore}
Suppose $ g \ge 0$, $0  \le h \le 1$, and $ h= 1$ on $\{ x: g(x) >0\} $. 
If $ \rho$ is a p-conic quasi-linear functional then $ \rho(ag+bh) =a\rho(g) + b\rho(h)$ for any $a, b \ge 0$.
If $ \rho$ is a quasi-linear functional, then $ \rho(ag+bh) =a\rho(g) + b\rho(h)$ for any $a, b \in \mathbb{R}$.
\item \label{rhoCcCon}
Suppose  $ \rho$ is a p-conic quasi-linear functional, $f, g \in C_c(X),  f, g \ge 0$. Then 
$| \rho(f) - \rho(g)| \le  \parallel \rho \parallel \parallel f-g \parallel$.
\end{enumerate} 
If $X$ is compact we also have:
\begin{enumerate}[label=(\roman*),ref=(\roman*)]
\item \label{addCcl}
If $ f \in C(X)$, and a functional $\rho$ is conic-linear on $A^+(f)$ 
(or on $A^-(f) $), then $ \rho(f + c) = \rho(f) + \rho(c)$. 
\item
If $\rho$ is conic-linear on $A^+(h)$ for each function $h \in C(X)$ (respectively, 
conic-linear on $A^-(h)$ for each function $h \in  C(X)$)  and $ \rho(1)  \in \mathbb{R}$, 
then $\rho$ is monotone. 
\end{enumerate}
\end{lemma}  

\begin{proof}
\begin{enumerate}[label=(\Roman*),ref=(\Roman*)]
\item
Since $0 \in A^+(f)$, $\rho$ is conic-linear on $A^+(f)$, and $| \rho(0) | < \infty$ we have  $\rho(0) = 0$. A similar argument works for $ A^-(f)$.
\item
Suppose $ fg=0, f \ge 0, g \le 0 $. Taking $h = f + g$, observe that $ f  = (id \vee 0) \circ h, \, g = (id \wedge 0) \circ h$, 
so   $f, g \in A^+(h)$. 
Then $\rho(f+g) = \rho(f)  + \rho(g)$.
\item 
Note that $ c \ge 0$ and $f= (id \wedge c) \circ (f+g), \ g=(0 \vee (id-c)) \circ (f+g)$.
\item
We may assume that $ g \neq 0$.
Since $ 0 \le g \le \parallel g \parallel h \le \parallel g \parallel$, by part \ref{levcoL} 
$$ \rho(ag + bh) = \rho(ag + \frac{b}{\parallel g \parallel} \parallel g \parallel h) = a \rho(g) + \frac{b}{\parallel g \parallel} \rho(\parallel g \parallel h)
= a \rho(g) + b \rho(h).$$
\item
Suppose  $f, g \in C_c(X), \, f,g \ge 0$.
Choose $ h \in C_c(X)$  such that $ h \ge 0,  h = 1$ on $ supp \, f \cup supp \, g$.
Since $f-g \le \parallel f-g \parallel \,  h  $, i.e. $ f \le g + \parallel f-g \parallel \,  h $, 
using Definition \ref{cqlf} and part \ref{levcoMore} we have:
\[ \rho(f) \le \rho(g + \parallel f-g \parallel \, h ) =  \rho(g) + \parallel f-g \parallel \, \rho(h). \] 
Similarly, $ \rho(g) \le \rho(f) + \parallel f-g \parallel \ \rho(h)$, so 
\begin{eqnarray}  \label{rhofgh}
| \rho(f ) -\rho(g) | \le  \parallel f-g \parallel \,   \rho(h)  \le  \parallel f-g \parallel  \ \parallel \rho \parallel.
\end{eqnarray}
\end{enumerate}
\begin{enumerate}[label=(\roman*),ref=(\roman*)]
\item
Since every constant $c \in A^+(f)$ (and also $c \in A^-(f)$) we 
immediately see that $ \rho(f + c) = \rho(f) + \rho(c)$. 
\item
Note that  $ \rho(c) \in \mathbb{R}$ for any constant $ c \ge 0$.   
Let $g \le f$. Choose a constant $c>0$ such that $ 0 \le g+c \le f+c$.  We have:
$$ \rho(g) + \rho(c) = \rho(g+c) \le \rho(f+ c) = \rho(f) + \rho(c), $$
so $ \rho(g) \le \rho(f)$.
\end{enumerate}
\end{proof}

\begin{proposition} \label{Rphointor}
Suppose $X$ is locally compact and $\rho$ is a functional on $C_0(X)$ that is positive-homogeneous and monotone on 
nonnegative functions. If $ \rho$ is real-valued, then $\| \rho \| < \infty$.
\end{proposition}

\begin{proof}
The statement can be obtained by adapting the argument from Lemma 2.3 in \cite{Alf:ReprTh}
(see Proposition 50 in \cite{Butler:QLFLC}).
\end{proof} 

\section{d-functionals} \label{SeDfnl}

In this section we shall define several functionals. The domains of these functionals vary, but the most common are 
$C(X), C_0(X), C_c(X), C_0^+(X)$. We shall not specify the domains in the definitions and results that hold for different domains,
but shall indicate them later when we use these functionals on specific collections of functions.

\begin{definition} \label{dfnl}
A functional $\rho$ with values in $[ -\infty, \infty]$ (assuming at most one of $\infty, - \infty$) and $| \rho(0) | < \infty$ 
is called a d-functional if   
on nonnegative functions it is positive-homogeneous, monotone, and orthogonally additive, i.e. for $f, g \in D(\rho)$
\begin{enumerate}[label=(d\arabic*),ref=(d\arabic*)]
\item 
$f \ge 0, \ a > 0  \Longrightarrow  \rho (a f) = a \rho(f)$, 
\item 
$0 \le  g \le f \Longrightarrow  \rho(g) \le \rho(f) $,
\item  \label{ortadd}
$f \cdot g = 0, f,g \ge 0  \Longrightarrow  \rho(f + g) = \rho(f) + \rho(g)$. 
\end{enumerate}
\end{definition}

\begin{remark} \label{sdPos}
It is easy to see that $\rho(0) = 0$, and so $\rho$ is positive. 
\end{remark}

\begin{definition} 
We say that a functional $\rho$ satisfies the c-level condition if 
 $$ f \le c, \, f=c \mbox{  on  } supp \, g  \Longrightarrow  \rho(f + g) = \rho(f)  + \rho(g), $$ 
where $ c $ is a constant, and $f, g \in D(\rho)$. 

If the domain of $\rho$ includes constants, we say that $\rho$ satisfies the constant condition if
for any function $g \in D(\rho)$ and any constant $c$ 
\begin{align} \label{constUsl1}
\rho(g+ c) = \rho(g) + \rho(c).
\end{align} 
\end{definition}

\begin{lemma} \label{FcAddDfnl}
Suppose the domain of $\rho$ includes constants. 
\begin{enumerate}[label=(\roman*),ref=(\roman*)] 
\item 
The c-level condition implies the constant condition. 
\item \label{DfnlCadit}
Suppose $\rho$ is a d-functional that satisfies the constant condition, $\rho(1) \in \mathbb{R}$, and $c \ge 0$. 
If $ f \le c,  f=c$ on $ supp \, g$, $g \ge 0$,  then $\rho(g-f) = \rho(g)  + \rho(-f)$.
\item \label{Clev=Cnst}
Suppose a functional $\rho$ is conic-linear on $A^+(h)$ for each function $h$ (for example, $\rho$ is a 
p-conic quasi-linear). If f $\rho$ satisfies the constant condition
then $\rho$ satisfies the c-level condition for all $g \ge0$.
\end{enumerate}
\end{lemma}

\begin{proof}
\begin{enumerate}[label=(\roman*),ref=(\roman*)]
\item
Take $f \equiv c$.
\item
Since $c-f \ge 0$ and $(c-f) g = 0$, we have:
$$\rho(c) + \rho(g-f)  = \rho(c + g - f) = \rho(c-f) + \rho(g) 
= \rho(c) + \rho(-f) + \rho(g),$$
so  $\rho(g-f) = \rho(g)  + \rho(-f)$.
\item
Suppose $\rho$ satisfies the constant condition (\ref{constUsl1}). Let $ f \le c, \ f=c$ on $ supp \, g, \, g \ge 0$.
Let $h=f-c$.  Then $h \le 0$ and $g \,  h =0$.  By part \ref{addPO}  of Lemma \ref{cqiConst}, $\rho(h+g) = \rho(h) + \rho(g)$.
Then
\begin{align*}
\rho(f+g) &= \rho(h+c+g) = \rho(h+g) + \rho(c) \\
&= \rho(h)  + \rho(g) + \rho(c) = \rho(h+c) + \rho(g) = \rho(f) + \rho(g).
\end{align*}
\end{enumerate}
\end{proof} 

\begin{remark} \label{consNorm}
If $\rho$ is a monotone, positive-homogeneous functional that satisfies the constant condition 
then for any functions $f,g \in D(\rho)$
\[ | \rho(f) - \rho(g) | \le \rho(\parallel f-g \parallel ) = \parallel f-g \parallel  \rho(1).  \]
One can show this by noticing that $ f \le g + \parallel f-g \parallel $ and using an argument similar to the one 
for part \ref{rhoCcCon} of Lemma \ref{cqiConst}.
\end{remark}

We modify condition \ref{ortadd} in Definition \ref{dfnl}  and obtain the definition of a c-functional:

\begin{definition} \label{cfnl}
A functional $\rho$ with values in $[ -\infty, \infty]$ (assuming at most one of $\infty, - \infty$) and $| \rho(0) | < \infty$ 
is called a c-functional if for $f, g \in D(\rho)$ 
\begin{enumerate}[label=(c\arabic*),ref=(c\arabic*)]
\item \label{cfl1} 
$f \ge 0, \ a > 0 \Longrightarrow  \rho (a f) = a \rho(f)$;  
\item 
$0 \le  g \le f \Longrightarrow  \rho(g) \le \rho(f) $;
\item  \label{ortadd1}
$f \cdot g = 0, f,g \ge 0 \mbox {   or   } f \ge 0, g \le 0  \Longrightarrow  \rho(f + g) = \rho(f) + \rho(g)$.
\end{enumerate}
\end{definition}
Note that for a c-functional $\rho$, $\ \rho(f) = \rho(f^+) + \rho(-f^-)$.

\begin{definition} \label{sfnl}
A real-valued functional $\rho$ is called an s-functional if
\begin{enumerate}[label=(s\arabic*),ref=(s\arabic*)]
\item
$f \ge 0, \ a \in \mathbb{R} \Longrightarrow  \rho (a f) = a \rho(f)$;  
\item \label{OrdPres} 
$0 \le g \le f \Longrightarrow  \rho(g) \le \rho(f) $;
\item 
$f \cdot g = 0, f,g \ge 0 \mbox {   or   } f \ge 0, g \le 0  \Longrightarrow  \rho(f + g) = \rho(f) + \rho(g)$.
\end{enumerate}
\end{definition}

\begin{lemma} \label{eqDefSfnl}
Suppose $\rho$ is an s-functional. Then 
\begin{enumerate}[label=(\roman*),ref=(\roman*)]
\item \label{par111}
$\rho(f) =  \rho(f^+) - \rho(f^-). $
\item
$ \rho(af) = a \rho(f)$ for any $f$ and any $a \in \mathbb{R}$.
\item
If $ g \le f$ then  $ \rho(g) \le \rho(f)$.
\item
If $f\cdot g =0$ then $ \rho(f+g) = \rho(f) + \rho(g)$. 
\end{enumerate}
\end{lemma}

\begin{proof}
\begin{enumerate}[label=(\roman*),ref=(\roman*)]
\item
$\rho(f) = \rho(f^+ -f^-) = \rho(f^+) + \rho(-f^-) =  \rho(f^+) - \rho(f^-). $
\item 
Easy to see from part \ref{par111}.
\item
If $ g \le f$ then $g^+ \le f^+, \ f^- \le g^-$ and 
$ \rho(g) = \rho(g^+) - \rho(g^-) \le \rho(f^+) - \rho(f^-) = \rho(f).$
\item
If $ f \cdot g=0$ then
$ (f+g)^{+} = f^{+} + g^{+}, \   (f+g)^{-} = f^{-} + g^{-}, $ and 
$ f^{+} \cdot g^{+} =0,  \ f^{-} \cdot g^{-} =0 .$ 
Then 
\begin{eqnarray*}
\rho(f+g) & =& \rho((f+g)^{+}) - \rho((f+g)^{-})  
= \rho(f^{+} + g^{+}) -\rho(f^{-} + g^{-})  \\
&=& \rho(f^{+}) + \rho(g^{+}) - \rho(f^{-}) - \rho(g^{-}) = \rho(f) + \rho(g). 
\end{eqnarray*}
\end{enumerate}
\end{proof}

Lemma \ref{eqDefSfnl} shows that Definition \ref{sfnl} is equivalent to the following:

\begin{definition} \label{sfnl2}
A functional $\rho$ is called an s-functional if
it is homogeneous, monotone, and orthogonally additive, i.e.
\begin{enumerate}[label=(sa\arabic*),ref=(sa\arabic*)]
\item \label{sfnl2us1}
$ \rho (a f) = a \rho(f)$  for any $a \in \mathbb{R}$; 
\item \label{OrdPres1} 
$g \le f \Longrightarrow  \rho(g) \le \rho(f) $;
\item 
$f \cdot g = 0  \Longrightarrow  \rho(f + g) = \rho(f) + \rho(g)$.
\end{enumerate}
\end{definition}

\begin{remark} \label{QLFisSFNL}
Let $X$  be locally compact, and let $\rho$ be a quasi-linear functional. By~\cite[Lemma 20(q2), Section 3 and Lemma 44, Section 5]{Butler:QLFLC},  
$\rho$ is an s-functional. 
\end{remark}

\begin{definition} \label{PhiFam}
Let $QI$ and $L$ denote, respectively, the families of all quasi-linear and linear functionals. 
By $\Phi^d, \Phi^c, \Phi^s, \Phi^+, \Phi^- $ we denote, respectively, the families of all d-functionals, c-functionals, s-functionals, 
p-conic quasi-linear functionals, and n-conic quasi-linear functionals.   
\end{definition}

\begin{remark} \label{PhiFam1}
We have: $ L \subseteq QI \subseteq \Phi^s \subseteq \Phi^c \subseteq \Phi^d.$
\end{remark}

\begin{proposition} \label{dfnSFNL} 
\begin{enumerate}[label=(\roman*),ref=(\roman*)]
\item  \label{sfkon}
$\rho$ is a s-functional iff 
it is a real-valued c-functional with a property that $\rho(-g) = - \rho(g)$ for every $g \le 0$ in the domain of $\rho$.
\item
If $\rho$ is a d-functional, $D(\rho)$ contains constants, $\rho(1) \in \mathbb{R}$, and the constants condition (\ref{constUsl1}) is satisfied, then 
$ \rho$  is monotone and  positive.
\end{enumerate}
\end{proposition}

\begin{proof}
\begin{enumerate}[label=(\roman*),ref=(\roman*)]
\item 
Suppose $ \rho$ is a c-functional with the property that $\rho(-g) = - \rho(g)$ for every $g \le 0$ in the domain of $\rho$.
We have: $\rho(f) = \rho(f^+) + \rho(-f^-)  = \rho(f^+) - \rho(f^-)$ 
for any function $f$ in the domain of $\rho$.
Then $\rho(-f) = - \rho(f) $ for any function $f$ in the domain of $\rho$, and  condition \ref{sfnl2us1} of Definition \ref{sfnl} follows.
Thus, $ \rho$ is an s-functional.
\item
Let $g \le f$.
It is enough to assume that  $ 0 \le g \le f$, for
choosing a positive constant $k$ such that $ g+ k \ge 0$ we see that $\rho(g) \le \rho(f) $ iff  
$\rho(g) + \rho(k) \le \rho(f) + \rho(k)$ iff $\rho(g+k) \le \rho(f+k)$. 
But $\rho$ is a d-functional,  so $\rho(g) \le \rho(f)$ for $ 0 \le g \le f$.
Since $ \rho(0) =0$, monotonicity of $\rho$ implies positivity. 
\end{enumerate}
\end{proof}

Now we will introduce the closely related r- and l- functionals. 

\begin{definition} \label{rfnl}
A functional $\rho$ with values in $[ -\infty, \infty]$ (assuming at most one of $\infty, - \infty$) and $| \rho(0) | < \infty$  
is called an r-functional if  for $f, g \in D(\rho)$
\begin{enumerate}[label=(r\arabic*),ref=(r\arabic*)]
\item 
$ a > 0   \Longrightarrow   \rho (a f) = a \rho(f)$; 
\item 
$ 0 \le g \le f   \Longrightarrow  \rho(g) \le \rho(f) $;
\item  \label{ortaddR}
If $f \cdot g = 0$ where $ f,g \ge 0$ or $ f \ge 0, g \le0$ then  $ \rho(f + g) = \rho(f) + \rho(g)$.
\end{enumerate}
If $D(\rho)$ contains constants then we also require for $ f \in D(\rho) $ and a constant $c$
$$ \rho(f + c) = \rho(f) + \rho(c). $$ 
\end{definition}

\begin{definition} \label{lfnl}
A functional $\rho$ with values in $[ -\infty, \infty]$ (assuming at most one of $\infty, - \infty$) and $| \rho(0) | < \infty$  
is called an l-functional if 
\begin{enumerate}[label=(l\arabic*),ref=(l\arabic*)]
\item 
$ a > 0   \Longrightarrow   \rho (a f) = a \rho(f)$;
\item 
$  g \le f  \le 0  \Longrightarrow  \rho(g) \le \rho(f) $;
\item  \label{ortaddL}
If $f \cdot g = 0$ where $ f,g \le 0$ or $ f \ge 0, g \le0$ then  $ \rho(f + g) = \rho(f) + \rho(g)$.
\end{enumerate}
If $D(\rho)$ contains constants then we also require for $ f \in D(\rho) $ and a constant $c$
$$ \rho(f + c) = \rho(f) + \rho(c).$$ 
\end{definition}

\begin{definition} 
By $\Phi^r, \Phi^l $ we denote, respectively, the families of all r-functionals, and l-functionals on $X$.  
\end{definition}

\begin{remark} \label{RLD}
Here are a few easy observations.
\begin{enumerate}[label=(\roman*),ref=(\roman*)]
\item \label{agreeP}
Suppose $\rho, \mathcal{R}$ are r-functionals with the same domain that contain constants and $ \rho(1) \in  \mathbb{R}$. Then $ \rho = \mathcal{R}$ iff 
$ \rho = \mathcal{R}$ on nonnegative functions.
Indeed,  
for an arbitrary function $f $  choose a constant $c \ge 0$ such that $f + c \ge 0$ and see that 
 $ \rho(f) + \rho(c) = \rho(f+c) = \mathcal{R}(f+c ) = \mathcal{R}(f) + \mathcal{R}(c) = \mathcal{R}(f) + \rho(c)$, i.e.  $ \rho(f) = \mathcal{R}(f)$.
\item \label{fsINfr}
We have: $ \Phi^s \subseteq \Phi^r \subseteq \Phi^c \subseteq \Phi^d$ and from Definition \ref{sfnl2} $ \Phi^s \subseteq \Phi^l$.   
\item \label{Pifnl}
Given a functional $\rho$, consider also the functional $\pi$ defined by $\pi(f) = - \rho(-f)$ for every $f \in D(\rho)$.  
Then $\pi$ is an l-functional iff
$ \rho$ is an r-functional. 
\end{enumerate}
\end{remark}

\begin{lemma} \label{+isR}
Each p-conic quasi-linear functional is an r-functional. Each n-conic quasi-linear functional is an l-functional. 
So $\Phi^+ \subseteq \Phi^r$ and $ \Phi^- \subseteq \Phi^l$.
\end{lemma}

\begin{proof}
Using part \ref{addPO} and part \ref{addCcl} of Lemma \ref{cqiConst} we see that 
a p-conic quasi-linear functional is an r-functional. 
The second statement follows from Remark \ref{PiA+} and part \ref{Pifnl} of Remark \ref{RLD}.
\end{proof}

\begin{lemma} \label{CcondRfl}
Let $X$ be compact and $ \rho$ be a functional on $C(X)$.  
\begin{enumerate}[label=(\Roman*),ref=(\Roman*)]
\item  \label{RLmon}
If $\rho$ is an r-functional with $ \rho(1) \in \mathbb{R}$  or an l-functional with $ \rho(-1) \in \mathbb{R}$  then $\rho$ is monotone.
\item \label{normforf}
If $\rho$ is an r-functional  with $ \rho(1) \in \mathbb{R}$ 
then $ | \rho(f) - \rho(g) | \le \parallel \rho \parallel  \, \parallel f-g \parallel = \rho(1) \parallel f-g \parallel.$
\item \label{clevel}
If $\rho$ is an r-functional with $ \rho(-1) \in \mathbb{R}$ then $\rho$ satisfies the c-level condition for any $ g \ge 0$.
\item \label{addnegf}
If $ \rho$ is an r-functional,  $f \cdot g = 0, \  f \ge 0 $ then  $ \rho(f + g) = \rho(f) + \rho(g)$.
Similarly, 
if $ \rho$ is an l-functional, $f \cdot g = 0, \  f \le 0 $ then  $ \rho(f + g) = \rho(f) + \rho(g)$.
\end{enumerate}
\end{lemma}

\begin{proof}
\begin{enumerate}[label=(\Roman*),ref=(\Roman*)]
\item
Suppose $ g \le f$. Choose a constant $c>0$ such that $0 \le g+ c \le f + c$. Then 
$ \rho(g) + \rho(c) = \rho(g+c) \le \rho(f+ c) = \rho(f) + \rho(c)$, 
which gives $ \rho(g) \le \rho(f).$ The monotonicity of an l-functional can be proved similarly.
\item
Use part \ref{RLmon} and Remark \ref{consNorm}. 
\item
Assume that $f \le c, f=c $ on $ supp \, g, g \ge 0$. Then 
$(f-c)g=0, f-c \le 0$, and 
$$ \rho(-c) + \rho(f+g) = \rho(f+g-c) = \rho(f-c) + \rho(g) = \rho(f)+ \rho(-c) + \rho(g) .$$ 
Thus, $\rho(f+g) = \rho(f) + \rho(g)$.
\item
Let $ \rho$ be an r-functional, $f \cdot g = 0, \,  f \ge 0 $. Then $(f + g^+)(-g^-) = 0, f \cdot g^+ =0$, and
by Definition \ref{rfnl}
\begin{align*}
 \rho(f+g) &= \rho(f + g^+ -g^-) = \rho(f+g^+) + \rho(-g^-) = \\
 & \rho(f) +\rho(g^+) + \rho(-g^-) = \rho(f) + \rho(g).
\end{align*}
The proof for an l-functional is similar.
\end{enumerate}
\end{proof}

\begin{remark}
Part \ref{addnegf} of Lemma \ref{CcondRfl} and part \ref{Pifnl} of Remark \ref{RLD}  were first observed for a compact space 
in~\cite[Propositions 17 and 18]{Svistula:DTM}.
\end{remark}

\begin{definition} \label{rgDef}
Let $\rho$ be a d-functional. Let $g \ge 0$. Define 
$\rho_g (f) = \rho(fg)$ for $f \in D(\rho)$. 
\end{definition}

\begin{lemma} \label{rg}
If $\rho$ is a  d-,c-, r-, l-,s- functional  or a linear functional, then so is $\rho_g$.
\end{lemma}

\begin{proof}
Easy to check.
\end{proof}

\section{Deficient topological measures from d-functionals} \label{SeGetDTM}

\begin{definition} \label{mrDfnl}
Let $X$ be locally compact, and let $\rho$ be a d-functional with 
$  C_c^+(X) \subseteq D(\rho) \subseteq C_b(X)$. 
Define a set function $\mu_{\rho} : \mathscr{O}(X) \cup \mathscr{C}(X) \rightarrow [0, \infty]$ as follows:
for an open set $U \subseteq X$ let 
\[ \mu_{\rho}(U) = \sup\{ \rho(f): \  f \in C_c(X), 0\le f \le 1, supp \, f \subseteq U  \}, \]
and for a closed set $F \subseteq X$ let
\[ \mu_{\rho}(F) = \inf \{ \mu_{\rho}(U): \  F \subseteq U,  U \in \mathscr{O}(X) \}.\]
\end{definition}

\noindent
Note that Definition \ref{mrDfnl} is consistent for clopen sets.

\begin{lemma} \label{PRmrDfnl}
For the set function $\mu_{\rho}$ from Definition \ref{mrDfnl} the following holds: 
\mbox{ } 
\begin{enumerate}[label=y\arabic*.,ref=y\arabic*]
\item \label{nongt1}
$\mu_{\rho} $ is nonnegative.
\item \label{opmon1}
$\mu_{\rho}$ is monotone.
\item \label{rhoK1}
Given an open set $U$, for any compact $K \subseteq U$ 
\[ \mu_{\rho}(U) = \sup \{ \rho(g):  1_K \le g \le 1, \ g \in C_c(X), \  supp \, g \subseteq U  \}. \]
\item \label{surho1}
For any $K \in \mathscr{K}(X)$ 
\[ \mu_{\rho}(K) = \inf \{ \rho(g): \   g \in C_c(X), g \ge 1_K \}. \]
\item \label{surho1a}
For any $K \in \mathscr{K}(X)$ 
\[ \mu_{\rho}(K) = \inf \{ \rho(g): \   g \in C_c(X), \, 1_K \le g \le 1 \}. \]
\item \label{VVbar1}
Given $ K \in \mathscr{K}(X)$, for any open $U$ such that $ K \subseteq U$
\[ \mu_{\rho}(K) = \inf \{ \mu_{\rho}(V): \ V \in \mathscr{O}(X), \ K \subseteq V \subseteq \overline V \subseteq U \}. \]
\item \label{innerrg1}
For any  $U \in \mathscr{O}(X)$
\[ \mu_{\rho}(U) = \sup \{ \mu_{\rho} (K): K \in \mathscr{K}(X), \ K \subseteq U \}. \]
\item \label{mropad1}
For any disjoint compact sets $K$ and $C$ 
\[ \mu_{\rho}(K \sqcup C) = \mu_{\rho}(K) + \mu_{\rho}(C). \]
\item \label{mropclad1}
Suppose $X$ is compact, $\rho(1)  \in \mathbb{R}$, and $\rho$ is an s-functional 
satisfying the constant condition (\ref{constUsl1}).
If $ K \subseteq U, \ K \in \mathscr{K}(X), \ U \in \mathscr{O}(X)$ then
\[ \mu_{\rho}(U) = \mu_{\rho}(K) + \mu_{\rho}(U \setminus K). \]
\end{enumerate}
\end{lemma}

\begin{proof}  
For part \ref{nongt1}, $\mu_{\rho}$ is nonnegative since $\rho$ is a positive functional by Remark \ref{sdPos}. Part \ref{opmon1} is easy to see.
Proofs for parts \ref{rhoK1} - \ref{mropad1} follow proofs of the corresponding parts 
of~\cite[Lemma 35, Sect. 4]{Butler:QLFLC}. We shall show part \ref{mropclad1}.

Let $ K \subseteq U, \ K  \in \mathscr{K}(X), \ U \in \mathscr{O}(X)$. 
First we shall show that 
\begin{align} \label{subadd6}
\mu_{\rho}(U \setminus K) +\mu_{\rho}(K) \ge  \mu_{\rho}(U).
\end{align}
If $\mu_{\rho}(K) = \infty$,  the inequality (\ref{subadd6}) trivially holds, so we assume that $\mu_{\rho}(K) < \infty.$ 
By Lemma \ref{easyLeLC} let $V \in \mathscr{O}(X) $ with compact closure be such that
\[ K \subseteq V \subseteq \overline V \subseteq U.\] 
For $\epsilon>0$  choose $W_1 \in \mathscr{O}(X)$ such that $K \subseteq W_1 \subseteq V $ and $\mu_{\rho}(W_1) < \mu_{\rho}(K) + \epsilon$. 
Also, there exists $W \in \mathscr{O}(X) $ with compact closure such that
 \[ K \subseteq W \subseteq \overline W \subseteq W_1 \subseteq V  \subseteq \overline V \subseteq U.\]
Choose an Urysohn function $ g \in C_c(X)$ such that $1_{\overline W} \le g \le 1, \  supp \,  g \subseteq W_1$. 
Then 
\[ \rho (g) \le \mu_{\rho}(W_1) < \mu_{\rho}(K) + \epsilon. \] 
First assume that $\mu_{\rho}(U) < \infty$.
By part \ref{rhoK1} choose $f \in C_c(X)$ such that $1_{\overline V} \le f \le 1,  \ supp \, f \subseteq U$, and
\[ \rho(f)  > \mu_{\rho}(U) - \epsilon.\]
Note that $0 \le f-g \le 1$, and, since $f-g = 0$ on $ \overline W$, we have $supp \, (f -g) \subseteq U \setminus K$. 
Since $ f=1$ on $supp \, g$ and $\rho$ is an s-functional,  by Remark \ref{PhiFam1} and part \ref{DfnlCadit} of Lemma \ref{FcAddDfnl} 
$$ \rho(g-f) = \rho(g) + \rho(-f) =  \rho(g) - \rho(f), $$
so 
\begin{align} \label {vazhNer}
\rho(f-g) = \rho(f) -\rho(g).
\end{align} 
Then we have:
$$ \mu_{\rho}(U \setminus K)   \ge \rho(f-g) = \rho(f) - \rho(g) \ge \mu_{\rho}(U) - \epsilon - \mu_{\rho}(K) - \epsilon,$$
which gives us inequality (\ref{subadd6}).
If $\mu_{\rho}(U) = \infty$, use instead of $f$ functions $f_n $ with $1_{\overline V} \le f_n \le 1,  \ supp \, f_n \subseteq U, \ \rho(f_n) \ge n$
in the above argument to show that $\mu_{\rho}(U \setminus K) = \infty$. Then inequality (\ref{subadd6}) holds.

Now we would like to show that 
\begin{align} \label{supadd6a} 
\mu_{\rho}(U) \ge  \mu_{\rho}(U\setminus K) + \mu_{\rho}(K).
\end{align}
By monotonicity of $\mu_{\rho}$ it is enough assume that $\mu_{\rho}(U \setminus K), \mu_{\rho}(K)  < \infty$.
Given $\epsilon>0$,  
choose $ g \in C_c(X),  \ 0 \le g \le 1$ such that $ C = supp \, g \subseteq U \setminus K$ and 
\[ \rho(g) > \mu_{\rho}(U \setminus K)  - \epsilon.\] 

Note that $K \subseteq U \setminus C$.
If $\mu_{\rho}(U \setminus C) = \infty,$ then $\mu_{\rho}(U) = \infty$, so (\ref{supadd6a}) holds. 
So assume that $\mu_{\rho}( U \setminus C) < \infty.$ 
By part \ref{rhoK1} choose $f \in C_c(X)$ such that 
$1_K \le f \le 1,  \ supp \, f \subseteq U \setminus C$, and $ \rho(f)  > \mu_{\rho}(U \setminus C) - \epsilon$.
Then 
\[ \rho(f)  > \mu_{\rho}(U \setminus C) - \epsilon \ge \mu_{\rho}(K) - \epsilon.\]
Since $fg=0, \, f,g \ge 0$,  we have
$\rho(f+g) = \rho(f) + \rho(g)$. Since $ f+g \in C_c(X)$ with  $supp \,(f+ g) \subseteq U$, we obtain:
\begin{align*}
\mu_{\rho}(U) \ge \rho(f+g)  = \rho(f) + \rho(g) \ge \mu_{\rho}(K)  + \mu_{\rho}(U \setminus K) - 2 \epsilon.
\end{align*}
Therefore,   $ \mu_{\rho}(U) \ge  \mu_{\rho}(U\setminus K) + \mu_{\rho}(K)$. 
\end{proof} 

\begin{remark} \label{subaddQL}
In the proof of part \ref{mropclad1} the only place where we 
need the fact that $X$ is compact and $D(\rho)$ contains constants
is when we use part \ref{DfnlCadit} of Lemma \ref{FcAddDfnl} to obtain formula (\ref{vazhNer}) in order to get 
inequality (\ref{subadd6}).
If $\rho$ is a quasi-linear functional on a locally compact space then formula (\ref{vazhNer})
holds by part \ref{levcoL} of Lemma \ref{cqiConst}, and we again obtain inequality (\ref{subadd6}). 
Our means of obtaining inequality (\ref{subadd6}) resembles one from~\cite[Theorem 3.9]{Alf:ReprTh}.
\end{remark} 

\begin{theorem} \label{rho2muD}
Suppose  $X$ is locally compact, $\rho$ is a d-functional with 
$  C_c(X) \subseteq D(\rho) \subseteq C_b(X)$, and  $\mu_{\rho}$ defined in Definition \ref{mrDfnl}.  
Then
\begin{enumerate}[label=(\roman*),ref=(\roman*)]
\item
$\mu_{\rho}$ is a deficient topological measure. 
\item
If $\rho$ is real-valued on  $ C_c(X)$, then $\mu_{\rho}$ is compact-finite.
\item
If $ \rho $ is bounded, then $\mu_{\rho}$ is finite.
\item
If $X$ is compact and $D(\rho) = C(X)$  then $\mu_{\rho}(X) = \rho(1)$.
\item \label{s1tm}
If the domain of $\rho$ includes constants, $\rho(1) \in \mathbb{R}$, and $\rho$ is an s-functional satisfying  
constant condition (\ref{constUsl1}) then  $\mu_{\rho}$ is a topological measure.
\item \label{QLdastTM}
If $ \rho$ is a quasi-linear functional then  $\mu_{\rho}$ is a topological measure.
\end{enumerate}
\end{theorem}

\begin{proof} 
\begin{enumerate}[label=(\roman*),ref=(\roman*)]
\item
Note that since  $\rho$ is not identically $\infty$, then neither is $\mu_{\rho}$. 
By part \ref{nongt1} of Lemma \ref{PRmrDfnl} $\mu_{\rho}$ is nonnegative.
Part \ref{mropad1} of Lemma \ref{PRmrDfnl} gives \ref{DTM1} of Definition \ref{DTM}. 
Definition \ref{mrDfnl} and part \ref{innerrg1} of Lemma \ref{PRmrDfnl} give regularity conditions  \ref{DTM2} and 
\ref{DTM3} of Definition \ref{DTM}. Thus, $\mu_{\rho}$ is a deficient topological measure.
\item
Follows from part \ref{surho1} of Lemma \ref{PRmrDfnl}.
\item
Evident from Definition \ref{mrDfnl}.
\item
See Definition \ref{mrDfnl}.
\item
Follows from part \ref{mropclad1} (or just inequality (\ref{subadd6})) of Lemma \ref{PRmrDfnl} and Theorem \ref{DTMtoTM}. 
\item
Follows form Remark \ref{subaddQL} and Theorem \ref{DTMtoTM}. 
\end{enumerate}
\end{proof} 

\section{Left and right measures} \label{lerim}

Given  a deficient topological measure and a bounded continuous function  we may consider four 
distribution functions.

\begin{definition} \label{4fns}
Let $\mu$ be a finite deficient topological measure on a locally compact space $X$. Let 
$f  \in C_b(X)$. Define the following nonnegative functions on $\mathbb{R}$:
\begin{align*}
L_1 (t) &=L_{1,\mu, f} (t) = \mu(f^{-1} ((-\infty, t) )), \\
L_2 (t) &=L_{2, \mu, f} (t) = \mu(f^{-1} ((-\infty, t] )), \\
R_1 (t) &= R_{1, \mu, f} (t) =  \mu(f^{-1} ((t, \infty) )), \\
R_2 (t) &=  R_{2,  \mu, f} (t) =\mu(f^{-1} ([t, \infty) )).
\end{align*}
\end{definition}

\begin{remark} \label{L1longform} 
For particular $\mu$ and $f$  to simplify notations we use $L_1, L_2, R_1, R_2$. 
When we need to emphasize the dependence on 
$\mu$ and $f$, we use notations $L_{1,\mu, f},   R_{1, \mu, f} $ and so on. 
When we need to use, say, $L_1$ as a function of $f$ we denote it by $L_{1, f}$. 
\end{remark} 

\begin{lemma} \label{4distrib}
Let $\mu$ be a finite deficient topological measure on a locally compact space $X,\ f  \in C_b(X)$. 
Let  nonnegative real-valued functions  $L_1, L_2, R_1, R_2$ be as in Definition \ref{4fns}. Then  
\begin{enumerate}[label=\Roman*.,ref=\Roman*]
\item \label{cha3}
Functions $L_1, L_2$ are non-decreasing; $R_1, R_2$ are non-increasing. 
If $f(X) \subseteq [a,b]$ then 
$$ L_1(a) = L_2(a^-) = 0; \ \ \ \ \  L_1(b^+) =  L_2(b^+)=  L_2(b) = \mu(X), $$
$$ R_1(a^-) =R_2(a^-) = R_2(a) = \mu(X); \ \ \ \ \   R_1(b) = R_2(b^+) = 0.$$ 
\item \label{cha2}
$L_1$ is left-continuous,  $R_1$ is right-continuous.
\item \label{cha4}
$L_1(t^-) = L_2(t^-)  = L_1(t)$ for any $t$.  
If $L_2$ is left-continuous at $t$ (in particular, continuous at $t$) then $L_1(t) =L_2(t)$. Similarly, 
$R_2(t^+) = R_1(t^+) = R_1(t)$ for any $t$, and if $R_2$ is right-continuous at $t$ then $R_1(t) = R_2(t)$. 
In particular,  the set of $t$ where $L_1(t) \neq  L_2(t) $ and  the set of $t$ where $R_1(t) \neq  R_2(t) $ are, at most, countable sets. 
\item \label{cha5}
$L_1(t) + R_1(t) \le \mu(X)$ for every $t$.
\item \label{cha6}
If $X$ is compact, the function $L_2$ is right-continuous, and $R_2$ is left-continuous.
If $X$ is locally compact, $ f \in C_0(X)$, then $R_2(t)$ is left-continuous at $a$ and any $t>0$, and  $L_2(t)$ is right-continuous at $b$ and any $t<0$.
In particular, $R_2$ is left-continuous at all $t$ except, possibly, $t \in E$ for some countable set $E \subseteq (-\infty, 0] \setminus \{ a\}$ and 
$L_2(t)$ is right-continuous at all $t$ except, possibly,  $t \in E_1$ for some countable set $E_1 \subseteq [0, \infty) \setminus \{b\}$.
\end{enumerate}
\end{lemma}

\begin{proof}
\begin{enumerate}[label=\Roman*.,ref=\Roman*]
\item \label{cha3p}
Easy to see.
\item \label{cha2p}
The sets $U_s = f^{-1} ((-\infty,s) )$ are open,  $U_s \nearrow U_t$ as $s \rightarrow t^{-}$, 
so by Lemma \ref{opaddDTM} $L_1$ is left-continuous. The argument for $R_1$ is similar.
\item \label{cha4p}
Let $s < t$. Then $f^{-1}((-\infty, s)) \subseteq f^{-1}((-\infty, s] ) \subseteq f^{-1}((-\infty, t))$, and so  
$L_1(t^-) \le L_2(t^-)  \le L_1(t) \le L_2(t)$.  By left-continuity of $L_1$ we have $L_1(t^-) = L_2(t^-)  = L_1(t)$;   
if  $L_2$ is left-continuous at $t$  
then $L_1(t) = L_2(t)$. Similarly for $R_1$ and $R_2$. 
\item \label{cha5p}
The sets $ f^{-1} ((-\infty, t) )$ and $ f^{-1} ((t, \infty) )$ are disjoint open sets, so 
from superadditivity of $\mu$ we see that $L_1(t) + R_1(t) \le \mu(X)$
for every $t$. 
\item   
If $X$ is compact, the sets $C_a = f^{-1} ([a, \infty) $ are compact. 
From  Lemma \ref{opaddDTM} it follows that $R_2$ is left-continuous. 
If $X$  is locally compact and $ f \in C_0(X)$ then the sets $K_a = f^{-1} ([a, \infty)), \ a>0 $ are compact. 
From  Lemma \ref{opaddDTM} it follows that $R_2$ is left-continuous at any $t >0$. The assertions about $L_2$ are proved similarly.
\end{enumerate}
\end{proof}

\begin{remark} \label{2integrals}
Let $ \mu$ be a finite deficient topological measure on a locally compact space $X$. Let $f \in C_b(X)$ with $ f(X) \subseteq [a,b]$.
By Theorem \ref{LebSt} and part \ref{cha3} of Lemma \ref{4distrib} the Riemann-Stieltjes integral
$$ \int_a^b id \, dL_1 = - \int_a^b L_1(t) dt + L_1(b^+) b =  - \int_a^b L_1(t) dt + b \mu(X).$$
Let $l$ be the Lebesque-Stieltjes measure associated with $L_1$, so $l$ is a regular Borel measure on $ \mathbb{R}$.
By part \ref{cha4} of Lemma \ref{4distrib} we see that 
$$  \int_a^b id \, dl = \int_a^b id  \, dL_1 =  - \int_a^b L_1(t) dt + b \mu(X) = - \int_a^b L_2  (t) dt + b \mu(X).  $$
Let $r$ be the Lebesque-Stieltjes measure associated with $-R_1$, a regular Borel measure on $ \mathbb{R}$.
We have: 
$$  \int_a^b id \,  dr = \int_a^b id \, d(-R_1) =   \int_a^b R_1 (t) dt + a \mu(X) =  \int_a^b R_2 (t) dt + a \mu(X).   $$
 \end{remark}

\begin{definition} \label{rlmery}
We call $l$ the left measure and $r$ the right measure.
When the right and left measures are equal, we set $m=r=l$.
\end{definition}  

\begin{remark} \label{rlmeryLong}
The measures $r$ and $l$ arise from functions  $R_1 = R_{1, \mu, f}$ and $L_1  =L_{1,\mu, f}$. 
We use notations $R_{1, f}, r_f, l_f$  when we need to emphasize the dependence of $R_1$ and 
measures $r, l$ on the function $f$. If we want to use measures 
$r$ and $l$ as functions of $f$ and $\mu$, we write $r_{f, \mu}, \ l_{f, \mu}$.

When $\mu$ is a topological measure, measure $m$ is equal to  $\mu_f$ in~\cite{Alf:ReprTh} and 
$m_f$ in~\cite{Butler:QLFLC}. See~\cite[Remark 28, Sect. 3]{Butler:QLFLC}.
\end{remark} 

\begin{theorem} \label{RLDTM}
Let $\mu$ be a finite deficient topological measure on a locally compact space $X$, and let 
$f  \in C_0(X)$. 
\begin{enumerate}[label=(\Roman*),ref=(\Roman*)]
\item \label{rints}
There are regular Borel measures $r$ and $l$ on $\mathbb{R}$  such that 
$ supp \, \, r \subseteq \overline{ f(X )},  supp \, \,  l \subseteq  \overline{ f(X )}, \  l(\mathbb{R}) = \mu(X), r(\mathbb{R}) = \mu(X)$, 
$$ r((t, \infty)) = \mu( f^{-1} ((t, \infty)))  \mbox{    for all   } t, $$
$$ r([t, \infty)) = \mu( f^{-1} ([t, \infty)))  \mbox{    for all   } t \notin E, $$
where $E \subseteq (-\infty, 0]$ is a countable set from part \ref{cha6} of Lemma \ref{4distrib}
$$ l((- \infty, t)) =   \mu( f^{-1} ((- \infty, t)))   \mbox{    for all   } t,$$
$$ l((- \infty, t]) =   \mu( f^{-1} ((- \infty, t]))  \mbox{    for all   } t \notin E_1,$$
where $E_1 \subseteq [0, \infty)$ is a countable set from part \ref{cha6} of Lemma \ref{4distrib}.
\item  \label{rints2}
For any open or closed set $A  \subseteq \mathbb{R}$ 
$$ \mu(f^{-1} (A)) \le l(A),\ \ \  \mu(f^{-1} (A)) \le r(A). $$
\end{enumerate}
\end{theorem} 

\begin{proof}
Let $ \overline{ f(X)} = [a,b]$.
\begin{enumerate}[label=(\Roman*),ref=(\Roman*)]
\item
By Lemma \ref{4distrib} $L_1$ is left-continuous, so for every $t$
$$ l((-\infty, t)) = L_1(t) = \mu(f^{-1} ((- \infty, t))).$$

Next, using Lemma \ref{opaddDTM}
$$
l(\mathbb{R}) = \lim_{n \rightarrow  \infty} r((- \infty, n)) =  \lim_{n \rightarrow \infty}  \mu( f^{-1}((- \infty, n)) ) = \mu(f^{-1} (\mathbb{R})) = \mu(X).
$$

If $t \notin E_1$, 
then $L_2$ is right-continuous at $t$. Since $L_2 = L_1$ outside of 
a countable set, $l(( -\infty, t]) = \lim_{s \rightarrow t^+} l(( -\infty, s)) =  \lim_{s \rightarrow t^+} L_1(s) = 
 \lim_{s \rightarrow t^+} L_2(s)  = L_2(t) = \mu(f^{-1}( ( -\infty, t]))$.

Since $L_1$ is constant on $(-\infty, a) $ and on $ (b, \infty)$, we see that 
$l((-\infty, a)) =l((b, \infty)) =0$. It follows that $supp \, l \subseteq [a,b] =  \overline{ f(X)}$.

The statements for $r$ can be proved similarly.
\item
Let $(a,b) \subseteq \mathbb{R}, b \notin E$. By the superadditivity of $\mu$ and part \ref{rints}
\begin{align*}
\mu(f^{-1} ((a, b)) & \le  \mu(f^{-1} ((a, \infty))) - \mu(f^{-1} ([b, \infty))) \\
&= r((a, \infty)) - r([b, \infty)) =   r((a, b)).
\end{align*}
For  $(a,b) $ with $b \in E$, choose $b_n \notin E $ such that $(a,b_n) \nearrow (a,b)$. 
Since by  Lemma \ref{opaddDTM}  $ \mu \circ f^{-1} $ and $r$ are both 
$\tau$-smooth on open sets, we have $ \mu(f^{-1} ((a,b))) \le r ((a,b))$ for any $(a,b)$. 
We see that  $ \mu(f^{-1} (J)) \le r (J)$ for any finite or infinite open interval $J$. 
Then the same inequality holds for any open set $W \subseteq \mathbb{R}$.

Now let $C \subseteq \mathbb{R}$ be closed. Then 
\begin{align*}
\mu(f^{-1} (C))  &\le \inf \{\mu(f^{-1} (W) : \ C \subseteq W, \, W \in \mathscr{O}(X) \} \\
& \le \inf \{r(W) : \ C \subseteq W, \, W \in \mathscr{O}(X) \}  =r(C).
\end{align*} 

The statements for the left measure $l$ can be proved in a similar way.
\end{enumerate}
\end{proof}

\begin{theorem} \label{RLravny}
Let $ \mu $ be a finite deficient topological measure on a locally compact space. Then $ r= l$  iff $ L_1(t) + R_1(t) = \mu(X)$ 
for a.e. $t$ with respect to the Lebesque measure.
\end{theorem}

\begin{proof}
Using Remark \ref{2integrals} and part \ref{cha5} of Lemma \ref{4distrib} we may note that 
\begin{align*}
r = l  & \Longrightarrow \int_a^b id \, dL_1  =  \int_a^b id \, d(-R_1)   \\
& \iff \int_a^b R_1 dx + a \mu(X)  = - \int_a^b L_1 dx + b \mu(X)   \\
& \iff \int_a^b (L_1  + R_1) dx = (b-a) \mu(X)  \\
& \iff L_1 + R_1 = \mu(X) \mbox{   a.e.},
\end{align*}
where a.e. is with respect to the Lebesque measure $ \lambda$.

Conversely, 
let $L_1(t) + R_1(t) = \mu(X)$  for $t \notin D$, where $\lambda(D) = 0$. We may assume that $D$ contains sets $E, E_1$ from
part \ref{cha6} of Lemma \ref{4distrib} and 
all points where $L_1 \neq L_2, R_1 \neq R_2$. 
If $[a, b] \subseteq \mathbb{R}$ and $ a,b \notin E$ then  by part \ref{rints} of Theorem \ref{RLDTM} we have:
\begin{eqnarray*}
l ([a,b]) = l((-\infty, b]) - l((-\infty, a))  = L_2(b) - L_1(a)  \\
= \mu(X) - R_1(b) - \mu(X) + R_2(a) =  r([a, \infty)) - r((b, \infty)) = r([a,b]).
\end{eqnarray*}
An arbitrary interval $(a,b)$ can be written as $ \bigcup_{n=1}^{\infty} [a_n, b_n]$, where intervals $ [a_n, b_n]$ are ordered by inclusion, 
and  $a_n, b_n \notin E$. It follows that  measures $l= r$ on $\mathbb{R}$.
\end{proof}

\begin{theorem} \label{RLforTM}
Let $ \mu$ be a finite topological measure on a locally compact space $X$. Let $f \in C_0(X)$. Then for the right and left measures 
$r, l$  we have $r=l$.
\end{theorem}

\begin{proof}
Let $f(X) \subseteq [a,b]$. By Theorem \ref{RLravny} it is enough to show that there is a countable set $D$ 
such that $ L_1(t) + R_1(t) = \mu(X)$  for all $t \in \mathbb{R} \setminus D$.
Let $D$ be the countable (by part \ref{cha4} of Lemma \ref{4distrib}) set consisting of $0$ and  all points where
$L_1 \neq L_2, R_1 \neq R_2$. 
If  $t >0$ then $f^{-1} ([t, \infty))$ is compact, and it follows from \ref{TM1} of Definition \ref{TMLC} that
$L_1(t) + R_2(t) = \mu(X)$.  $ R_1(t) = R_2(t)$ for $t \in (0, \infty) \setminus D$,  
so $ L_1(t) + R_1(t) = L_1(t) + R_2(t) = \mu(X).$
Similarly, for all  $t \in (-\infty, 0) \setminus D$ we have $ L_1(t) + R_1(t) = L_2(t) + R_1(t) = \mu(X).$
\end{proof}

\begin{theorem} \label{mf2t}
Let $\mu$ be a finite topological measure on a locally compact space $X$, and let regular Borel measure $m  = r = l$.
\begin{enumerate}[label=(\Roman*),ref=(\Roman*)]
\item \label{mfDTMc0}
If $f  \in C_0(X)$, 
then $m (A) = \mu ( f^{-1}(A))$ for any open set $ A \subseteq \mathbb{R}$  and any closed set  $A  \subseteq   \mathbb{R} \setminus \{0\}$. 
\item
If $X$ is compact, $f \in C(X)$  then 
$m (A) = \mu ( f^{-1}(A))$ for any open  or closed set $ A \subseteq \mathbb{R}$. 
\end{enumerate} 
\end{theorem}

\begin{proof}
\begin{enumerate}[label=(\Roman*),ref=(\Roman*)]
\item \label{mfDTMcs}
First let $ (a,b) \in \mathbb{R}, \, b \neq 0 $.
Note that in  
\[  f^{-1}( (a, \infty)) = f^{-1}( (a,b) ) \sqcup f^{-1}( \{b\} ) \sqcup f^{-1}( (b, \infty) ), \]
all the sets are open except for the middle set on the right hand side, which is compact since $f \in C_0(X)$. 
Applying $\mu$ we obtain 
\begin{align} \label{FFF}
R_1(a) = \mu( f^{-1}( (a,b) ))  + \mu(f^{-1}( \{ b \})) + \mu(f^{-1}( (b, \infty)))
\end{align}
Since $ f^{-1}( (b, \infty)) \sqcup f^{-1}( \{ b \} ) \subseteq  f^{-1}( (t, \infty) ) $ for any $t < b$,
by superadditivity (see Lemma \ref{opaddDTM}) we have:
\[  \mu(f^{-1}( \{\ b \} )) + \mu(f^{-1}( (b, \infty))) \le \mu( f^{-1}( (t, \infty) )) = R_1(t). \]
Thus,  from (\ref{FFF}) we see that $R_1(a) \le   \mu(f^{-1}( (a,b) )) + R_1(t). $
As  $ t \to b^-$ we have:
\[ m((a,b)) = R_1(a) - R_1(b^-) \le \mu(f^{-1}( (a,b) )).\]
Together with part \ref{rints2} of Theorem \ref{RLDTM}  we obtain $ m((a,b)) = \mu(f^{-1}( (a,b) )) $ for any interval $(a,b), \, b \neq 0$. 
An interval $(a,0) = \cup_{n=1}^{\infty} (a, -\frac 1n)$. Since both $\mu$ and $m$ are $\tau-$ smooth and
additive on open sets (see Lemma \ref{opaddDTM}), the result  holds for any  
finite open interval in $\mathbb{R}$, and then for any open set in $\mathbb{R}$. 
Below in \ref{mfDTM2b} we shall prove that $\mu(f^{-1}(C))  = m(C)$ 
for closed sets in $\mathbb{R} \setminus \{0\}$.

\item \label{mfDTMcomp}
The set $f^{-1}([d, \infty) )$ is compact for every $d$, and the argument as in part \ref{mfDTMcs} shows that 
$m (W) = \mu ( f^{-1}(W))$ for every open set $W \subseteq \mathbb{R}$. 

Now let $C \subseteq \mathbb{R}$ be closed. Choose $W \in \mathscr{O}(X)$ such that $ C \subseteq W$.
Since $f^{-1}(C)$ is compact and $\mu$ is a topological measure, 
\begin{eqnarray*}
m(W) -   \mu(f^{-1}(C)) = \mu(f^{-1}(W)) - \mu(f^{-1}(C)) \\
= \mu(f^{-1}(W \setminus C)) = m (W \setminus C).
\end{eqnarray*}
Thus  $\mu(f^{-1}(C)) = m (C)$.

\item \label{mfDTM2b}
Now we shall finish the proof of part \ref{mfDTMc0}.
Let  $C \in \mathbb{R} \setminus \{0\}$.
Set $C_1 = C \cap (0, \infty)$ and $C_2 = C \cap (-\infty, 0)$. We have $C = C_1 \sqcup C_2$, and $b = \inf C_1 >0$. 
Since $f^{-1} (C_1) \subseteq f^{-1} ([b, \infty)) $, the  set $f^{-1} (C_1)$ is compact.  
An argument similar to the one in part \ref{mfDTMcomp}
shows that  $\mu(f^{-1}(C_1))  = m (C_1)$. Similarly, $\mu(f^{-1}(C_2))  = m (C_2)$, 
and so by finite additivity of $\mu$ and $m $ on compact sets $\mu(f^{-1}(C))  = m (C)$.  
\end{enumerate}
\end{proof}

\begin{lemma} \label{abPtMa}
Let $\mu$ be a finite deficient topological measure on a locally compact space $X$.
\begin{enumerate}[label=\Roman*.,ref=\Roman*]
\item
If $\mu $ is a simple deficient topological measure, then measures $r$ and $l$ are point masses, $l=\delta_a, \ r= \delta_b, \, b \le a$, 
where 
$$ a= \inf \{ t : L_1(t) = 1\} = \sup\{ s: L_1(s) = 0 \},$$ 
$$ b = \inf\{ t : R_1(t) = 0\} = \sup\{ s: R_1(s) = 1\}.$$ 
\item \label{p2cpfn}
If $X$ is compact and  $f =c$ is a constant function, then the measure $m = \mu(X) \delta_c$, where $\delta_c$ is a point mass at $c$.
\end{enumerate}
\end{lemma}

\begin{proof}
\begin{enumerate}[label=\Roman*.,ref=\Roman*]
\item
Since $\mu$ is simple, the non-decreasing function $L_1$ assumes only two values, and has single discontinuity at 
$a= \inf \{ t : L_1(t) = 1\} = \sup\{ s: L_1(s) = 0 \}$.
Since $l(\{a\}) = L_1(a^+) - L_1(a^-) = 1$, we see that $l = \delta_a$.
Similarly, $ r= \delta_b$, where $ b = \inf\{ t : R_1(t) = 0 \} = \sup\{ s: R_1(s) = 1\}$.

If $a < b$ then $L_1(t)  + R_1(t) = 2$ on interval $(a,b)$, which contradicts part \ref{cha5} of Lemma \ref{4distrib}. Thus, $b \le a$.
\item
We have $R_1(t) = \mu(X) $  for every $t < c$, and  $R_1(t) =0 $ for every $t \ge c$. 
Then $ m ( \{c\}) = r ( \{c\})  =  R_1(c^-) - R_1(c^+) = \mu(X)$.
\end{enumerate}
\end{proof}

\begin{example}  \label{baless}
Let $X = \mathbb{R}, \ D = [0,1]$ and $\mu$ be a simple deficient topological measure as in~\cite[Example 48, Sect. 6]{Butler:DTMLC}, i.e. 
$\mu(A) = 1$ if $ D \subseteq A$ and $\mu(A) = 0$ otherwise,  where $A \in \mathscr{O}(X) \cup \mathscr{K}(X)$.
Consider the following $f \in C_0(X): f(0) = 1, f(t) = 0$ for $ t \in (- \infty, -1] \cup [1, \infty)$, and $f$ is linear on $[-1,0]$ and $[0, 1]$.
Note that $D \subseteq f^{-1} ((- \infty, t))$ iff $t >1$, and  $D \subseteq f^{-1} ((t, \infty))$ iff $t <0$. Thus, 
\begin{eqnarray*}
 L_1 (t) & = &
  \left\{
  \begin{array}{rl}
  1 & \mbox{ if } t >1 \\
  0 & \mbox{ if } t \le 1 \\
  \end{array}
  \right.
\end{eqnarray*}
and  
\begin{eqnarray*}
 R_1 (t) & = &
  \left\{
  \begin{array}{rl}
  1 & \mbox{ if } t < 0 \\
  0 & \mbox{ if } t \ge 0 \\
  \end{array}
  \right.
\end{eqnarray*} 
So $l = \delta_1$ and $r=\delta_0$.
\end{example}

\begin{example} \label{baeq}
Let $X = \mathbb{R}$, the family $\mathcal{E} = \{D= [1,2]\}, \lambda_{0}  = \delta_{3/2}$. Let $\mu = \lambda^{+}$ as in~\cite[Example 49, Sect. 6]{Butler:DTMLC}.
Consider the following $f \in C_0(X): f(0) = 2, f(t) = 1 $ for $t \in [1,2], f(t) = 0$ for $ t \in (- \infty, -1] \cup [3, \infty)$, 
and $f$ is linear on $[-1,0], [0, 1]$, and $[2,3]$.
Note that $D \subseteq f^{-1} ((- \infty, t))$ iff $t >1$, and  $D \subseteq f^{-1} ((t, \infty))$ iff $t <1$. It follows that 
\begin{eqnarray*}
 L_1 (t) =\lambda^{+}(f^{-1} ((- \infty, t))) & = &
  \left\{
  \begin{array}{rl}
  1 & \mbox{ if } t >1 \\
  0 & \mbox{ if } t \le 1 \\
  \end{array}
  \right.
\end{eqnarray*}
and  
\begin{eqnarray*}
 R_1 (t) = \lambda^{+}(f^{-1} ((t, \infty)) ) & = &
  \left\{
  \begin{array}{rl}
  1 & \mbox{ if } t < 1 \\
  0 & \mbox{ if } t \ge 1 \\
  \end{array}
  \right.
\end{eqnarray*} 
Thus, for measures $l,r$ we have $ l = r = \delta_1$.
\end{example}

\begin{remark} \label{suppRL} 
In Theorem \ref{RLDTM} it is stated that  $ supp \, l, supp \, r \subseteq f(X)$. 
In Example \ref{baless}  and Example \ref{baeq}  $supp \, l $ and $ supp \, r $ are properly contained in $f(X)$. 
On the other hand, 
from part \ref{p2cpfn}  of Lemma \ref{abPtMa}  we see that it is also possible to have $ supp \, l  = supp \, r = f(X)$. 
\end{remark} 

\begin{remark} \label{abHow}
From  part \ref{cha5} of Lemma \ref{4distrib} we know that $L_1(t) + R_1(t) \le \mu(X)$. 
Although  $L_1(t) + R_1(t) = \mu(X)$ a.e. when $ \mu$ is a finite topological measure 
(see Theorems \ref{RLravny} and \ref{RLforTM}),
for deficient topological measures we may have both situations: in Example \ref{baeq}  $L_1(t) + R_1(t) = \mu(X)$ a.e., but in 
Example \ref{baless}  we have  $L_1(t) + R_1(t) = 0 < \mu(X) $ for $ t \in  (0,1)$.

In Lemma \ref{abPtMa} we have $b \le a$.  Example \ref{baless} and 
Example \ref{baeq} show that both situations when $b<a$ and $b=a$ are possible. 
These examples also show that when $\mu$ is a deficient topological measure, 
we can have both situations for measures $l$ and $r$ induced by $\mu$ and a given function $f$: when $l=r$ and when $l \neq r$.
\end{remark}

\section{Functionals from deficient topological measures} \label{SeGetFnls}

When $\mu$ is a finite deficient topological measure (not a topological measure) the measures $r, l$ are not equal in general, and 
we consider two different integrals:
$$  \int_a^b id \, dl $$ 
and    
$$  \int_a^b id  \, dr. $$   

\begin{definition} \label{DTM2rho}
Let $\mu$ be a finite deficient topological measure on a locally compact space $X$, and let measures 
$r=r_{f, \mu},l = l_{f, \mu}, m=m_{f, \mu}$ be 
as in Definition \ref{rlmery}, Remark \ref{2integrals}, and Remark \ref{rlmeryLong}.
Define  the following functionals on $ C_b(X)$:
$$ \mathcal{R} (f) = \mathcal{R}_{\mu} (f) =  \int _{\mathbb{R}}  id \, dr,$$
$$ \mathcal{L} (f) = \mathcal{L}_{\mu} (f) =  \int _{\mathbb{R}}  id \, dl,$$
and 
$$ \rho (f) = \rho_{\mu}(f) = \int _{\mathbb{R}}  id \, dm.$$
\end{definition}

\begin{remark}  \label{RHOforms} 
By Theorem \ref{RLDTM} $ supp \,  r , \, supp \,  l \subseteq \overline{f(X)}$,  so for any $[a,b]$ containing $f(X)$ 
$$ \mathcal{R} (f) = \int _{\mathbb{R}}  id \,  dr =  \int _{[a,b]}  id \, dr, \ \ \  \mathcal{L} (f) =   \int _{[a,b]}  id \, dl, \ \ \ \rho(f) =   \int _{[a,b]}  id \, dm. $$

With functions $L_1, L_2, R_1, R_2$ as in  Definition \ref{4fns} by Remark \ref{2integrals}  we have: 
\begin{align} \label{lfform}
 \mathcal{L} (f) &= \int _{\mathbb{R}}  id \,  dl  = \int_a^b id \, dl = - \int_a^b L_1 (t)  dt + b \mu(X) = - \int_a^b L_2 (t) dt + b \mu(X).  
\end{align} 
\begin{align} \label{rfform}
\mathcal{R} (f) & = \int _{\mathbb{R}}  id \,  dr = \int_a^b id \, dr  =   \int_a^b R_1 (t) dt + a \mu(X)  =  \int_a^b R_2 (t) dt + a \mu(X).  
\end{align}
If $ f=c$ is a constant function then 
\begin{align} \label{rfc}
\mathcal{R} (f)  = \mathcal{L}(f) = c \mu(X). 
\end{align}
If $f$ is nonnegative with $f(X) \subseteq [0,b]$ we have:
\begin{align} \label{rfformp}
 \mathcal{R} (f) = \int_0^b  id \, dr  =   \int_0^b R_1 (t) dt =   \int_0^b R_2 (t) dt.
\end{align}
When $r=l$ (in particular, when $\mu$ is a topological measure) and $f(X) \subseteq [0,b]$   we have 
 $$ \rho(f) = \int_0^b  id \, dm =   \int_0^b R_1 (t) dt =   \int_0^b R_2 (t) dt. $$
Similarly,  if $f$ is nonpositive with $f(X) \subseteq [a,0]$ we have:
\begin{align} \label{lfformn}  
\mathcal{L} (f) =  \int_a^0 id \, dl  =   - \int_a^0 L_1 (t) dt = -  \int_a^0 L_2  (t) dt.
\end{align}
When $r=l$ (in particular, when $\mu$ is a topological measure) and $f(X) \subseteq [a,0]$   we have 
 $$ \rho (f) = \int_a^0  id \, dm =  -  \int_a^0 L_1 (t) dt =  - \int_a^0 L_2 (t) dt. $$
\end{remark}

\begin{remark}
Note that when $\mu$ is a topological measure, we obtain familiar formulas.
See, for example,~\cite[Proposition 3.7]{Alf:ReprTh} and~\cite[Remark 43, Section 5]{Butler:QLFLC}. 
These results were, in turn, generalizations of results first obtained by J. F. Aarnes for compact spaces in~\cite{Aarnes:TheFirstPaper}. 
For example, when $X$ is compact and $\mu(X) = 1$, 
formula (\ref{rfform}) gives~\cite[formula (3.3)]{Aarnes:TheFirstPaper}.
\end{remark}

\begin{remark} \label{RLconnd}
We have the connection between $\mathcal{R}$ and $\mathcal{L}$ (which is the same as noted in~\cite[p. 739]{Svistula:DTM}).
We use  notations as in Remark \ref{rlmeryLong}. 
Observe that $R_{1,f} (-t) = L_{1,-f} (t)$. Thus, $l_{-f} = r_f \circ T^{-1}$, where $T(t) = -t$ for $ t \in \mathbb{R}$. 
Then $ \int_{\mathbb{R}} id \, dl_{-f} = -  \int_{\mathbb{R}} id \, dr_{f}$, 
i.e. 
\begin{align} \label{rlrr}
\mathcal{L}(-f) = - \mathcal{R}(f).
\end{align}
We may prove results for $\mathcal{R}$ and obtain similar results for $\mathcal{L}$ by analogy 
(as we did, for example, in Theorem \ref{RLDTM}) or 
using relation (\ref{rlrr}).  
\end{remark}

\begin{definition} \label{QIwrtDTM}
Let $ \mathcal{R}$ be the functional as in Definition \ref{DTM2rho}. We call the functional $\mathcal{R}$ a quasi-integral (with 
respect to a deficient topological measure $ \mu$) and write:
\begin{align*}
\int_X f \, d\mu = \mathcal{R}(f) = \mathcal{R}_{\mu} (f) =  \int _{\mathbb{R}}  id \, dr.
\end{align*}
\end{definition} 

\begin{remark}
If $\mu$ is a topological measure on $X$, by Definitions \ref{DTM2rho} and \ref{rlmery} we obtain exactly 
the quasi-integral in~\cite[Definition 27, Section 3]{Butler:QLFLC}.
\end{remark}

\begin{lemma} \label{LRsvva}
Let $\mathcal{L}, \mathcal{R}$ be functionals as in Definition \ref{DTM2rho}.
\begin{enumerate}[label=(\roman*),ref=(\roman*)]
\item \label{LRortad} 
$\mathcal{R}$ is orthogonally additive on nonnegative functions, and $\mathcal{L}$ is orthogonally additive on nonpositive functions.
\item \label{LROO} 
 $\mathcal{R}(0) = 0$ and $\mathcal{L} (0) = 0$.
\item \label{LRposh}
$\mathcal{L}, \mathcal{R}$ are positive-homogeneous functionals.
\item \label{LRmnt}
$\mathcal{L}, \mathcal{R}$ are monotone.
In particular, $\mathcal{L}, \mathcal{R}$ are positive. 
\end{enumerate}
\end{lemma} 

\begin{proof}
We use notations as indicated in Remark \ref{rlmeryLong}.
\begin{enumerate}[label=(\roman*),ref=(\roman*)]
\item
Let $f,g \ge 0$ and $f g =0$. Say, $f(X), g(X) \subseteq [0,b]$. 
For any $t >0$ observe that  $(f+g)^{-1} ((t, \infty)) = f^{-1} ((t, \infty))  \sqcup g^{-1} ((t, \infty))$, 
so by additivity of a deficient topological measure
on disjoint open sets we immediately obtain $R_{1, f+g} (t) =  R_{1, f} (t) + R_{1, g} (t)$. 
Since $(f+g)(X) \subseteq [0,b]$, from (\ref{rfformp}) we have $\mathcal{R}(f+g) = \mathcal{R}(f) + \mathcal{R}(g)$. 
Thus, $\mathcal{R}$ is orthogonally additive on nonnegative functions. Then 
orthogonal additivity of $\mathcal{L}$ on nonpositive functions follows from (\ref{rlrr}).
\item
Follows from part \ref{LRortad} or from (\ref{rfformp}) and (\ref{lfformn}). 
\item
If $c=0$ then from \ref{LROO} we see that $\mathcal{R}(cf) = \mathcal{L}(cf) = 0$. 
Let $c >0$. Since $L_{1, cf} (t) = L_{1,f} (t/c) $, from (\ref{lfform}) we see that $\mathcal{L}(cf) = c \mathcal{L}(f)$.
One can show that $\mathcal{R}$ is also positive-homogeneous in a similar way using (\ref{rfform}) or using positive-homogenuity 
of $\mathcal{L}$ together with formula (\ref{rlrr}). 
\item
Suppose that $f \le g$. Choose an interval $[a,b]$ which contains both $f(X)$ and $g(X)$. 
Since $L_{1, f} \ge L_{1,g}$ and $R_{1, f} \le R_{1,g}$, from  (\ref{lfform}) and  (\ref{rfform}) we see that 
$\mathcal{L} (f) \le \mathcal{L}(g) $ and $\mathcal{R}(f) \le \mathcal{R}(g)$. 
\end{enumerate}
\end{proof}

\begin{lemma} \label{RnaA+}
If $ h =  \phi \circ f \in  A^+(f) $ and $[a,b]$ is any interval containing $f(X)$  then 
$$\mathcal{R}(h) = \int_{[a,b]} \phi \, dr_f. $$
Similarly, if $ h =  \phi \circ f \in  A^-(f) $ then 
$$ \mathcal{L}(h) = \int_{[a,b]}  \phi \, dl_f.$$
\end{lemma}

\begin{proof}
Let $r_f, r_h$ be the right measures for functions $f$ and $h$ as in Theorem \ref{RLDTM},  $r_f$ is supported on $\overline{f(X)} \subseteq [a,b]$.
Since $\phi$ is nondecreasing, 
for any interval $(t, \infty)$ we have $\phi^{-1}((t, \infty)) = (s, \infty)$, where $s = \min\{ y \in [a,b]: \phi(y) = t \} $. 
(This is similar to~\cite[Proposition 13(1)]{Svistula:DTM}.)
Then 
\begin{align*}
r_h ((t,\infty)) &= \mu(h^{-1}((t, \infty))) = \mu(f^{-1}(\phi^{-1}((t, \infty)))) \\
& = \mu(f^{-1}((s, \infty)))  = r_f ((s, \infty)) = (r_f \circ \phi^{-1} )((t, \infty)).
\end{align*}
Thus,  $r_h = r_f \circ \phi^{-1}$ are equal as measures. 
Using formula (\ref{rfform}) we have:
$$ \mathcal{R}(h) = \int_\mathbb{R} id \, dr_h =\int_{[a,b]}  \phi \, dr_f.$$
The formula  $ \mathcal{L}(h) = \int_{[a,b]} \phi \, dl_f$ can be proved in a similar way.
\end{proof} 

\begin{lemma} \label{polinA+}
The functional $\mathcal{R}$ is conic-linear on each cone $ A^+(f)$, and the functional $\mathcal{L}$ is conic-linear on each cone $ A^-(f)$. 
\end{lemma}

\begin{proof}
Suppose $h, g \in  A^+(f), h = \phi \circ f, g = \psi \circ f, \overline{f(X)} \subseteq [a,b]$. 
Applying Lemma \ref{RnaA+} we have:
\begin{align*} 
\mathcal{R}(h+g) &= \mathcal{R}( \phi \circ f+  \psi \circ f) = \mathcal{R}(( \phi +\psi) \circ f) = \int_{[a,b]}  (\phi +\psi) \, dr_f \\
&=  \int_{[a,b]}   \phi  \, dr_f  +  \int_{[a,b]}  \psi \, dr_f= \mathcal{R}(h) + \mathcal{R}(g).
\end{align*}
Since by Lemma \ref{LRsvva} $\mathcal{R}$ is also positive-homogeneous, we see that $\mathcal{R}$ is conic-linear on $ A^+(f)$ for each $f$.  
The statements for $\mathcal{L}$ can be proved similarly.
\end{proof}

\begin{theorem} \label{RisA+}
\begin{enumerate}[label=(\roman*),ref=(\roman*)]
\item
The functional $\mathcal{R}$ is a p-conic quasi-linear functional, and the functional $\mathcal{L}$ is an n-conic quasi-linear functional.  
\item \label{RfisR}
The functional $\mathcal{R}$ is an r-functional , and $\mathcal{L}$ is an l-functional. 
\end{enumerate}
\end{theorem}

\begin{proof}
\begin{enumerate}[label=(\roman*),ref=(\roman*)]
\item
Follows from Lemma \ref{LRsvva} and Lemma \ref{polinA+}.
\item
Apply Lemma \ref{+isR}.
\end{enumerate}
\end{proof}

\begin{remark} \label{An:conicLin} 
When $\mu$ is a topological measure, the functional $ \rho = \rho_\mu$ is a quasi-linear functional 
(see~\cite[Theorem 30, Section 3]{Butler:QLFLC}), so $ \rho$ is linear on each singly generated subalgebra.
Theorem \ref{RisA+} gives an analog of this for the case when $\mu$ is a deficient topological measure: 
if $\mu$ is a deficient topological measure, 
then the functional $\mathcal{R}$ obtained from $\mu$ is p-conic linear, so, in particular,  it 
is conic-linear on the cone 
$  A^+(f)$ for each $f$. 
\end{remark}

The next lemma shows properties that relate $\mathcal{R}$ and $\mu$.

\begin{lemma} \label{MUrhoMU}
Let $\mu$ be a finite deficient topological measure, and $\mathcal{R}, \mathcal{L} $ be functionals on $ C_0(X)$obtained from $\mu$ 
as in Definition \ref{DTM2rho}. 
\begin{enumerate}[label=z\arabic*.,ref=z\arabic*]
\item \label{MUrhoMU1}
If $ U \in \mathscr{O}(X) $ and $ f \in C(X) $ is such that  $ supp \, f \subseteq U , \ 0 \le f \le 1$ then $\mathcal{R} (f) \le \mu(U)$.
\item \label{MUrhoMU2}
If $C \in \mathscr{C}(X)$ and $ f $ is such that $ 0 \le f \le 1, \ f =1$ on $C$, then 
$\mathcal{R} (f) \ge \mu(C)$.
\item \label{MUrhoMU2a}
For any $f \in C_0(X)$ 
\[ \mu(X)  \cdot \inf_{x \in X} f(x)  \le \mathcal{R}(f)  \le \mu(X) \cdot \sup _{x \in X} f(x).  \]
Similarly,
\[ \mu(X)  \cdot \inf_{x \in X} f(x)  \le \mathcal{L}(f)  \le \mu(X) \cdot \sup _{x \in X} f(x).  \]
Hence, $ \parallel \mathcal{R} \parallel, \parallel \mathcal{L} \parallel \le \mu(X) $.
\item \label{MUrhoMU3}
If $ U \in \mathscr{O}(X) $  then 
$ \mu(U) = \sup\{ \mathcal{R} (f): \  f \in C_c(X), 0\le f \le 1, supp \, f \subseteq U  \}. $
\item \label{MUrhoMU4}
If $K  \in \mathscr{K}(X)$  then 
$ \mu(K) = \inf \{ \mathcal{R} (g): \   g \in C_c(X), g \ge 1_K \}. $
\item \label{MUrhoMU5}
If $K  \in \mathscr{K}(X)$  then 
\begin{align*}
\mu(K) &= \inf \{ \mathcal{R} (g): \   g \in C_c(X), g \ge 1_K , 0 \le g \le 1\} \\
           &=  \inf \{ \mathcal{R} (g): \   g \in C_0(X), g \ge 1_K \}  \\
           &= \inf \{ \mathcal{R} (g): \   g \in C_0(X), g \ge 1_K, 0 \le g \le 1 \}.
\end{align*}            
\item \label{MUrhoMU6}
$ \parallel \mathcal{R} \parallel = \mu(X) = \parallel \mathcal{L} \parallel $, so 
$ | \mathcal{R}(f) | \le  \parallel \mathcal{R} \parallel \parallel f \parallel$ and 
$ | \mathcal{L}(f) | \le  \parallel \mathcal{L} \parallel \parallel f \parallel$.
\item \label{MUrhoMUdo}
If $f, g \in C_c(X), \, f,g \ge 0, \, supp \, f, supp \, g \subseteq K$ where $K$ is compact, then 
$$ | \mathcal{R}(f) - \mathcal{R}(g) | \le \parallel f - g \parallel \, \mu(K).$$ 
\end{enumerate}
If $X$ is compact  we also have:
\begin{enumerate}[label=(\roman*),ref=(\roman*)]
\item \label{MUrhoMU7}
If $c$ is a constant then  $\mathcal{R} (c) = c \mu(X)$ and, hence, $\mathcal{R} (f + c) = \mathcal{R} (f) + c \mu(X)$.
\item \label{MUrhoMU8}
For any functions $f,g \in D(\rho)$ 
\[ | \mathcal{R}(f) - \mathcal{R}(g) | \le \mathcal{R}(1) \parallel f-g \parallel = \mu(X) \parallel f-g \parallel.\]
\end{enumerate}
\end{lemma}

\begin{proof}
By Theorem \ref{RLDTM}, if $ 0 \le f \le 1$ then  $supp \, r \subseteq  [0,1]$.
\begin{enumerate}[label=z\arabic*.,ref=z\arabic*]
\item \label{MUrhoMU1P}
Using formula (\ref{rfform}) and part \ref{rints} of Theorem \ref{RLDTM}, we have:
$\mathcal{R} (f) = \int_{\mathbb{R}} id  \, dr \le  1 \cdot r(\{t: t >0\}) = \mu(f^{-1} (0, \infty)) \le \mu(U).$
\item \label{MUrhoMU2P}
Using part \ref{rints} of Theorem \ref{RLDTM} we have:
$$
\mathcal{R} (f) = \int_{\mathbb{R}} id  \, dr \ge  1 \cdot r (\{t: t =1\}) =\mu(f^{-1} ([1, \infty)) \ge  \mu(C). 
$$
\item \label{MUrhoMU2aP} 
Let $a =  \inf_{x \in X} f(x), \ b=  \sup _{x \in X} f(x)$.  It is enough to assume $ a \notin E$ 
(see part \ref{rints} of Theorem \ref{RLDTM}),
for otherwise we may
take $a_n \nearrow a, a_n \notin E$.  
Then
\begin{align*}
a \mu(X) &= a \mu(f^{-1}([a, \infty))) = a r([a, \infty)) = a r([a,b]) \\
& \le \int_a^b  id  \, dr = \mathcal{R}(f) \le b r([a, \infty)) = b \mu(X).
\end{align*}
Because of formula (\ref{rlrr}), the statement for $\mathcal{L}(f)$ also holds.
\item \label{MUrhoMU3P}
By part \ref{MUrhoMU1},   $ \sup\{ \mathcal{R} (f): \  f \in C_c(X), 0\le f \le 1, supp \, f \subseteq U  \} \le \mu(U)$. 
For $ \epsilon >0$ choose $ K \in \mathscr{K}(X) $ such that $ \mu(U) - \mu(K) < \epsilon$.
Pick $ f \in C_c(X)$ such that $ 0 \le f \le 1, \ supp \, f \subseteq U, \ f=1$ on $K$. 
Then  by part \ref{MUrhoMU2} 
$ \mathcal{R} (f) \ge  \mu(K) \ge \mu(U) - \epsilon$. 
It follows that  $ \sup\{ \mathcal{R} (f): \  f \in C_c(X), 0\le f \le 1, supp \, f \subseteq U  \} = \mu(U)$. 
\item \label{MUrhoMU4P}
Take any $U \in \mathscr{O}(X)$ containing $ K $.   
Taking an Urysohn function $f$ such that $ f=1$ on $K$, 
$ supp \, f \subseteq U$ and using part \ref{MUrhoMU1} 
we see that
 $\inf \{ \mathcal{R} (g): \ g \in C_c(X),  g \ge 1_K \} \le \mathcal{R}(f)   \le \mu_{\rho}(U)$. Taking the infimum over all 
open sets containing $K$ we have: 
\[ \inf \{ \mathcal{R} (g): \ g \in C_c(X), g \ge 1_K  \}  \le \mu_{\rho}(K). \]
To prove the opposite inequality, take any $ g \in C_c(X)$ such that $g \ge 1_K$.
Let $0< \delta < 1$.
Let 
$$U = \{ x:  \ g(x) > 1-\delta \}.$$
Then $U$ is open and $ K \subseteq U$.
Consider function $ h = \inf \{g,  1-\delta \} $, so $ h \in  C_c(X)$.  
Since $0 \le h \le g$ we have $\mathcal{R} (h) \le \mathcal{R}(g).$ 
Because $\displaystyle{\frac{h}{1- \delta}}  = 1$ on $U$,
for any function $f \in C_c(X), \ 0 \le f \le 1, supp \, f \subseteq U$ we have
$ f \le \displaystyle{\frac{h}{1- \delta}}$ and so by parts \ref{LRposh} and \ref{LRmnt} of Lemma \ref{LRsvva}
$$  \mathcal{R} (f) \le \mathcal{R} \left( \frac{h}{1-\delta} \right) = \frac{\mathcal{R}(h)}{1-\delta} \le \frac{\mathcal{R}(g)}{1-\delta} \ . $$
Then 
\begin{align*}
\mu (K) &\le   \mu (U) = \sup \left\{  \mathcal{R} (f) : f \in  C_c(X), \, 0 \le f \le 1_U, \,  supp \, f \subseteq U  \right\}  \\ 
 &\le \frac{\mathcal{R}(g)}{1-\delta} \ .
\end{align*} 
Thus, for any  $ g \in C_c(X)$ such that $g \ge 1_K$ and any $0 < \delta < 1$ 
 \[ (1 - \delta) \mu(K) \le \mathcal{R}(g). \]
Therefore,
$ \mu (K) \le \inf \{ \mathcal{R}(g): \ g \in C_b(X), g \ge 1_K  \}.$ 
\item \label{MUrhoMU5P}
In the argument for part \ref{MUrhoMU4}  we may use, respectively, 
$g \in C_c(X), g \ge 1_K , 0 \le g \le 1$ or $  g \in C_0(X), g \ge 1_K $ or  $ g \in C_0(X), g \ge 1_K, 0 \le g \le 1$.
\item \label{MUrhoMU6P}
By part \ref{MUrhoMU3} we see that $ \mu(X) \le \parallel \mathcal{R} \parallel$.
Using also part \ref{MUrhoMU2a} we have $ \parallel \mathcal{R} \parallel = \mu(X)$, and 
by formula (\ref{rlrr}) also $\parallel \mathcal{L} \parallel = \mu(X)$.
\item
It is enough to consider $K= supp \, f \cup supp \, g$. For any function $ h \in C_c(X)$ such that $ h \ge 0, h=1$ on $K$ as in 
formula (\ref{rhofgh}) in the proof of part \ref{rhoCcCon} of Lemma \ref{cqiConst} we have:
$$ | \mathcal{R}(f) - \mathcal{R}(g)| \le \parallel f - g \parallel \, \mathcal{R}(h).$$ 
Taking the infimum over functions $h$, by part \ref{MUrhoMU4} we obtain the assertion.
\end{enumerate}
\begin{enumerate}[label=(\roman*),ref=(\roman*)]
\item \label{MUrhoMU7P}
By formula (\ref{rfc}), $ \mathcal{R}(c) = c \mu(X)$, and the rest of the statement follows from Theorem \ref{RisA+}.
\item \label{MUrhoMU8P}
Since  $\mathcal{R}$ is monotone, an r-functional, and $ \mu(X) = \parallel \mathcal{R} \parallel  = \mathcal{R}(1)$ the statement follows from 
Remark \ref{consNorm}.
\end{enumerate}
\end{proof}

\begin{remark}
The proof of part  \ref{MUrhoMU4} is similar to the one for~\cite[Lemma 35(p4), Sect. 4]{Butler:QLFLC}, which, in turn, 
follows a proof from~\cite{Grubb:Lectures}.
\end{remark} 

\begin{remark}
Part \ref{MUrhoMUdo} of Lemma \ref{MUrhoMU} means that on $ C_0^+(X)$ the functional $ \mathcal{R}$ is continuous with respect to 
topology of uniform convergence on compacta.
\end{remark}

\begin{definition} \label{nfWayDef}
Let $f$ be a continuous function on $X$ on a locally compact space $X$. Consider
the $\sigma$-algebra of subsets of $X$  
$$ \Sigma_f = f^{-1} ( \mathscr{B} (\mathbb{R})), $$
where $ \mathscr{B}(\mathbb{R})$ are the Borel subsets of $\mathbb{R}$.
Let $r, l, m $ be measures on $\mathbb{R}$ as in Theorem \ref{RLDTM} .
On $ \Sigma_f $ define measure $n_r, n_l, n $ by setting for each $ E \in  \Sigma_f $
$$ n_r(E)  = r(f(E)), \ \ \ n_l(E) =  l(f(E)), \ \ \ n(E) =  m (f(E)) $$ 
\end{definition}

\begin{remark} \label{nfWay}
Definition \ref{nfWayDef} leads to another way to represent  functionals $\mathcal{R}$ and  ${\mathcal{L}}$. 
Let $ \Sigma_f $ and measures $ n_r, \ n_l$ be as in Definition \ref{nfWayDef}.
Then for any set $D \in \mathscr{B} (\mathbb{R})$ we have:
\[ n_r \circ f^{-1} (D) =  r (f(f^{-1} (D))) = r (D),\ \ \   n_l \circ f^{-1} (D) = l(D), \]
i.e. 
\begin{align} \label{nrnlSokr}
n_r \circ f^{-1} = r, \ \  n_l \circ f^{-1} = l   \mbox{   and   } n \circ f^{-1} = m \mbox{   on   }   \mathscr{B} (\mathbb{R}).
\end{align}
By formula (\ref{rfform})
\[  \mathcal{R} (f) =\int_{\mathbb{R}}id \, dr=  \int_{\mathbb{R}} id \  \, d(n_r \circ f^{-1}) = \int_X f \,  \, dn_r,  \]
i.e. 
\[ \mathcal{R}(f) = \int_X f \,  \, dn_r. \]
Similarly, 
\[ \mathcal{L}(f) = \int_X f \,  \, dn_l. \]
\end{remark} 

\begin{lemma} \label{nfMU}
Let $\mu$ be a finite deficient topological measure on a locally compact space $X$, $ f \in C_0(X)$, 
and $n_r, \ n_l$ be measures defined in Definition \ref{nfWayDef}.
\begin{enumerate}[label=(\roman*),ref=(\roman*)]
\item \label{nfMU1} $ n_r ( f^{-1} ((t, \infty)) = \mu( f^{-1} ((t, \infty)) \mbox{    for all   } t, $
$$ n_r ( f^{-1} ([t, \infty)) = \mu( f^{-1} ([t, \infty)) \mbox{    for all   } t \notin E, $$
where $E \subseteq (-\infty, 0]$ is a countable set from part \ref{cha6} of Lemma \ref{4distrib};
$$ n_l ( f^{-1} ((- \infty, t))  =   \mu( f^{-1} ((- \infty, t))   \mbox{    for all   } t,$$
$$ n_l ( f^{-1} ((- \infty, t]) =   \mu( f^{-1} ((- \infty, t])  \mbox{    for all   } t \notin E_1,$$
where $E_1 \subseteq [0, \infty)$ is a countable set from part \ref{cha6} of Lemma \ref{4distrib}.
\item
For any open or closed set $A  \subseteq \mathbb{R}$ 
$$ \mu(f^{-1} (A)) \le n_r (f^{-1} (A)), \ \  \mu(f^{-1} (A)) \le n_l (f^{-1} (A)). $$
\item \label{nfMU2}
If $\mu$ is a finite topological measure and $f  \in C_0(X)$, 
then also $n(f^{-1} (A)) = \mu ( f^{-1}(A)$ for any open or closed set  $A  \subseteq \mathbb{R} \setminus  \{0\}. $
\item \label{nfMU3}
If $\mu$ is a finite topological measure and $X$ is compact, then 
$n (f^{-1}(A)) = \mu ( f^{-1}(A))$ for any open or closed set $ A \subseteq \mathbb{R}$.
\end{enumerate}
\end{lemma}

\begin{proof}
Follows from (\ref{nrnlSokr}), Theorem \ref{RLDTM}, and Theorem \ref{mf2t}.
\end{proof}

\begin{remark} 
Definition \ref{nfWayDef} and Remark \ref{nfWay}  were first observed for the case of the compact space $X$ by M. Svistula, 
see~\cite[(3.4)]{Svistula:DTM}.
\end{remark}

\section{Representation Theorems for deficient topological measures} \label{SeReprT}

\begin{theorem} [Representation theorem]  \label{ReprThLC}
Let $\mu$ be a finite deficient topological measure on a locally compact space $X$.
\begin{enumerate}[label=(\roman*),ref=(\roman*)]
\item \label{RTdtmLC}
Then there exists a unique 
p-conic  quasi-linear functional  $\mathcal{R}$ on $C_0^+(X)$ of finite norm such that $\mu = \mu_{\mathcal{R}}$ 
and $ \parallel \mathcal{R} \parallel = \mu(X)$.
In fact, $ \mathcal{R} = \mathcal{R}_{\mu}$. 
\item \label{RTdtmC}
If $X$ is compact,  the unique functional $ \mathcal{R}$ can be taken to be a real-valued r-functional on $C(X)$.  
\end{enumerate}
Here $\mathcal{R}_{\mu}$ and $\mu_{\mathcal{R}}$ are as in Definition \ref{DTM2rho} and Definition \ref{mrDfnl}.
\end{theorem}

\begin{proof}
\begin{enumerate}[label=(\roman*),ref=(\roman*)]
\item
Let $\mu$ be a finite deficient topological measure on a locally compact space $X$, and 
let $\mathcal{R} = \mathcal{R}_{\mu}$ be a functional on $C_0^+(X)$ obtained from $ \mu$ as in Definition \ref{DTM2rho}. 
Note that  $ \mathcal{R} $ is a d-functional by 
Theorem \ref{RisA+} and Remark \ref{RLD}. Then $\mu_{\mathcal{R}}$ defined as in Definition \ref{mrDfnl} from $\mathcal{R}$ is a 
deficient topological measure by Theorem  \ref{rho2muD}.  
From part \ref{surho1} of Lemma \ref{PRmrDfnl}  and part  \ref{MUrhoMU4} of Lemma \ref{MUrhoMU} 
we see that $\mu(K) = \mu_{\mathcal{R}}(K)$ 
for every compact $K$. Thus, $\mu = \mu_{\mathcal{R}}$.  
By part  \ref{MUrhoMU6} of Lemma \ref{MUrhoMU}  $ \mu(X) = \parallel \mathcal{R} \parallel.$

Now we shall show the uniqueness. Let 
$\eta$ be another p-conic  quasi-linear functional on  $C_0^+(X)$ of finite norm such that $\mu_{\eta} = \mu_{\mathcal{R}}$.
Then for some $ a \ge 1$
\[ \parallel \eta \parallel \le a \mu(X) .\]
Let $f \in C_0^+(X)$. Both $\rho $ and $ \eta$ are positive-homogeneous, so we may assume that $\overline{ f(X)} = [0,1]$.  
For  $n \in \mathbb{N}$  let  $t_i = i / n, \ i=0, \ldots, n$. 
For $ i=1, \ldots, n$ consider
functions $\phi_i$ defined as follows: 
\begin{eqnarray*}
  \phi_i (t) & = &
  \left\{
  \begin{array}{rl}
  0 & \mbox{ if }  t \le t_{i-1}\\
  t - t_{i-1} & \mbox{ if }  t_{i-1}  < t <  t_i \\
  \frac 1n  & \mbox{ if }  t \ge  t_i \,  .
  \end{array}
  \right.
\end{eqnarray*} 
With functions $f_i = \phi_i \circ f$ we have $f = \sum_{i=1}^n f_i$.
Since each $\phi_i$ is non-decreasing, and $\phi_i(0) =0$, each $f_i \in A^+(f)$.
Since  $\mathcal{R}$ and $\eta$ are both p-conic quasi-linear functionals, we have 
\begin{align} \label{razbrho}
 \mathcal{R} (f) = \mathcal{R}( \sum_{i=1}^n f_i) = \sum_{i=1}^n \mathcal{R}(f_i),\ \ \ \eta(f) =  \sum_{i=1}^n \eta(f_i). 
\end{align}
By part \ref{MUrhoMU1} of Lemma \ref{MUrhoMU} $ \mathcal{R} (nf_1) \le \mu(X)$, 
so $0 \le \mathcal{R}(f_1) \le  \frac1n \mu(X) \le \frac{a }{n} \mu(X)$.  
We have $\  0 \le \eta(nf_1) \le \parallel \eta \parallel  \le a \mu(X)$, 
i.e. $ 0 \le \eta(f_1) \le \frac{a}{n} \mu(X)$.
Thus, 
\begin{eqnarray} \label{estf1}
| \mathcal{R}(f_1) - \eta(f_1) | \le \frac{a \mu(X)}{n} \ . 
\end{eqnarray}
For each $ i =2, \ldots, n$ let $K_i = f^{-1}([t_i, \infty)),  V_i = f^{-1}((t_{i-1}, \infty))$. 
Choose an open set $U_i \subseteq V_i$ such that $\mu(U_i) - \mu(K_i) < \frac 1n$ and then pick 
an Urysohn function $g_i \in C_c(X)$ such that $0 \le g_i \le \frac 1n, \ g_i = \frac 1n$ on $K_i$ and $supp \, g_i \subseteq U_i \subseteq V_i$.
Since $ng_i = 1$ on $K$ and $\mu = \mu_{\mathcal{R}} = \mu_{\eta}$, by part \ref{surho1} of Lemma \ref{PRmrDfnl} 
and Definition \ref{mrDfnl}
 $ \mu(K_i) \le \mathcal{R}(n g_i), \eta (n g_i) \le \mu(U_i) $, and so  $ | \mathcal{R}(ng_i) - \eta(ng_i)| \le \mu(U_i) - \mu(K_i)  \le  \frac 1n$. 
Then
\begin{eqnarray} \label{vsne}
 | \mathcal{R}(g_i) - \eta(g_i)| \le \frac{ 1}{n^2}.
\end{eqnarray}
Let $g= \sum_{i=2}^n g_i, \, h= \sum_{i=2}^n f_i $. 
Since $supp \, (g_3 + \ldots + g_n) \subseteq K_2$ and $ g_2 =1$ on $K_2$, by part \ref{levcoL} of Lemma \ref{cqiConst}
$\rho(g_2 + g_3 + \ldots g_n) = \rho(g_2) + \rho(g_3 + \ldots g_n)$. Similarly for $\eta$. By induction
\begin{eqnarray} \label{Retasum}
 \mathcal{R}(g) = \sum_{i=2}^n \mathcal{R}(g_i), \ \ \ \ \ \eta(g) = \sum_{i=2}^n \eta(g_i). 
\end{eqnarray}
Then
\begin{eqnarray} \label{estg}
| \mathcal{R}(g) - \eta(g) | \le \sum_{i=2}^n | \mathcal{R}(g_i) - \eta(g_i) | \le \frac{n-1}{n^2} < \frac 1n.
\end{eqnarray}  
Note that $\parallel g-h \parallel \le \frac 1n$, so by part \ref{rhoCcCon} of Lemma \ref{cqiConst}
\begin{eqnarray} \label{estgh}
| \mathcal{R}(g) - \mathcal{R}(h) | \le \frac{ \mu(X)}{n}, \ \ \ \ | \eta(g) - \eta(h) | \le \frac{ a \mu(X)}{n}.
\end{eqnarray}
Using (\ref{estf1}), (\ref{estgh}), and (\ref{estg}) we obtain:
\begin{align*}
| \mathcal{R}(f) - \eta(f) | &\le | \mathcal{R}(f_1) - \eta(f_1) | + | \mathcal{R}(h) - \eta(h) | \\
 &\le | \mathcal{R}(f_1) - \eta(f_1)| + | \mathcal{R}(h) - \mathcal{R}(g)| \\
 &+ | \mathcal{R}(g) - \eta (g) | + |\eta(g) - \eta(h)|  \\
 & \le \frac{ a \mu(X)}{n} +\frac{ \mu(X)}{n}+ \frac 1n+ \frac{  a \mu(X)}{n} \le  \frac 1n (3 a \mu(X) + 1).
\end{align*}
Thus, $\mathcal{R} = \eta$.

\item
Now let $X$ be compact. 
We shall show that the proof for part \ref{RTdtmLC} still applies, although the reasoning for some estimates is different. 
We define the functional $ \mathcal{R} = \mathcal{R}_{\mu} $ on $C(X)$. 
It is an r-functional by Theorem \ref{RisA+}. 
Then $ |  \mathcal{R}(1) |  = \parallel \mathcal{R} \parallel =\mu(X) < \infty$, 
and $ |\mathcal{R}(-1)| \le \parallel \mathcal{R} \parallel < \infty$. By monotonicity, 
$\mathcal{R}$ is real-valued.
Let $ \eta $ be another real-valued (hence, bounded by Proposition \ref{Rphointor}) 
r-functional on $C(X)$ such that $\mu = \mu_{\mathcal{R}} = \mu_{\eta}$.  
To show the uniqueness, by Remark \ref{RLD} it is enough to show that   $\rho(f) = \eta(f) $ for $ f \ge 0$. 
For the functions $f_i$  from the proof for part \ref{RTdtmLC} 
note the following: $supp \, (f_2 + \ldots +  f_n) \subseteq  f^{-1}([t_1, \infty)), \ \ \ f_1 = \frac1n$ on $  f^{-1}([t_1, \infty))$, thus
by part \ref{clevel} of Lemma \ref{CcondRfl} we may apply the c-level condition to obtain 
$$ \mathcal{R} (f_1 + f_2 + \ldots +f_n)  = \mathcal{R} (f_1) + \mathcal{R} (f_2 + \ldots + f_n).$$ 
Then by induction we may show that formula (\ref{razbrho}) holds for $\mathcal{R}$ and, similarly, for $ \eta$.
In the same manner,  by part \ref{clevel} of Lemma \ref{CcondRfl} and induction  we show that
formula (\ref{Retasum}) holds. 
Note that  (\ref{estf1}), (\ref{vsne}), and (\ref{estg}) hold as in the proof of part \ref{RTdtmLC}.
Estimations (\ref{estgh}) are valid by part \ref{normforf} of Lemma \ref{CcondRfl}. 
Now as in the end of the proof for part \ref{RTdtmLC}, 
we show that $\mathcal{R} = \eta$.  
\end{enumerate}
The proof is complete now.  
\end{proof}

\begin{remark}
Our inequality (\ref{vsne}) is inspired by a similar estimate in the proof of~\cite[Theorem 9]{Svistula:DTM}.
\end{remark}

\begin{definition}
Let $ \mathbf{L}, \mathbf{QI}, \mathbf{\Phi^s}, \mathbf{\Phi^+} ,\mathbf{\Phi^-}, \mathbf{\Phi^r}, \mathbf{\Phi^l}$ 
represent subfamilies of bounded functionals 
from $L, QI, \Phi^s, \Phi^+, \Phi^-, \Phi^r, \Phi^l$
respectively. 
We may indicate in parenthesis the domain of functionals. 
For example, $\mathbf{\Phi^+}  =  \mathbf{\Phi^+} (C(X))$ is the collection of all bounded p-conic quasi-linear functionals on $C(X)$. 
\end{definition}

\begin{definition}
Let $\mathbf{DTM}(X),  \mathbf{TM}  (X), \mathbf{M} (X)$ represent, 
respectively,  subfamilies of finite set functions from $DTM (X), TM (X), M (X)$.
\end{definition}

\begin{corollary}  \label{ReprThLC-}
Let $X$ be locally compact. 
\begin{enumerate}[label=(\roman*),ref=(\roman*)]
\item \label{bijeG}
There is a bijection $\Gamma: \mathbf{\Phi^+} (C_0^+(X)) \longrightarrow \mathbf{DTM}(X)$ given 
by $\Gamma(\mathcal{R}) = \mu_{\mathcal{R}}$. 
The inverse bijection $\Pi = \Gamma^{-1}$ is given by  $\Pi : \mathbf{DTM}(X) \longrightarrow \mathbf{\Phi^+}  (C_0^+(X)) $  
where $\Pi(\mu) = \mathcal{R}_{\mu}$.  
Here  $\mu_{\mathcal{R}}$  and $\mathcal{R}_{\mu}$ are according  to Definition \ref{mrDfnl} and Definition \ref{DTM2rho}.
\item \label{RTdtmLC-}
There is a bijection $\Delta:  \mathbf{DTM}(X) \longrightarrow \mathbf{\Phi^-}(C_0^-(X))$ given by $\Delta(\mu) = \mathcal{L}_{\mu}$, 
where $ \mathcal{L}_{\mu}$ is according to Definition \ref{DTM2rho}.
\item \label{RTdtmCMP}
If $X$ is compact, there is a bijection between $ \mathbf{DTM}(X)$ and $ \mathbf{\Phi^r}(C(X))$ given by $\mu \mapsto \mathcal{R}_{\mu}$,  
and a bijection between $ \mathbf{DTM}(X)$ and  $ \mathbf{\Phi^l}(C(X))$. 
\end{enumerate}
\end{corollary}

\begin{proof}
\begin{enumerate}[label=(\roman*),ref=(\roman*)]
\item
Follows from Theorem \ref{ReprThLC}.
\item
By Remark \ref{PiA+} there is a bijection between $ \mathbf{\Phi^+} (C_0^+(X))$ and  $\mathbf{\Phi^-}(C_0^-(X))$, 
so we obtain bijection $ \Delta$. 
By formula (\ref{rlrr}) $\,  -\mathcal{R}_\mu (-f) = \mathcal{L}_{\mu} (f)$. 
\item 
By Theorem \ref{ReprThLC} there is a bijection between $ \mathbf{DTM}(X)$ and $ \mathbf{\Phi^r}(C(X))$ 
given by $\mu \mapsto \mathcal{R}_{\mu}$.
By Remark \ref{RLD}  there is a bijection between  $ \mathbf{\Phi^r}(C(X))$ and  $ \mathbf{\Phi^l}(C(X))$. 
\end{enumerate}
\end{proof} 

\begin{corollary} \label{S=pn}
Suppose  $ \rho$ is a finite c-functional on $C_0(X)$, and $\rho  = \mathcal{R}_{\mu}$ on $ C_0^+(X)$, 
where $\mathcal{R}_{\mu}$ is a functional on $ C_0(X)$ 
for some finite deficient topological measure $ \mu$ as in Definition \ref{DTM2rho}. Then 
$\rho$ is an s-functional iff $\mathcal{L}_{\mu}(g) = \rho(g)$ for all $ g \in C_0^-(X)$, where $\mathcal{L}_{\mu}$ is a functional 
on $ C_0(X)$ as in Definition \ref{DTM2rho}.
\end{corollary}

\begin{proof}
If $\rho(g) = \mathcal{L}_{\mu}(g)$ for all $ g \in C_0^-(X)$ then  $ \rho(g) = \mathcal{L}_{\mu}(g) = - \mathcal{R}_{\mu}(-g) = -\rho(-g)$, 
i.e. $ - \rho(g) = \rho(-g)$. 
Since $ \mu$ is finite, $\rho$ is real-valued on $C_0^+(X)$, so $\rho$ is real-valued.  
Then $ \rho$ is an s-functional by part \ref{sfkon} of Proposition \ref{dfnSFNL}. The other direction can be proved similarly.
\end{proof} 

\begin{theorem} \label{DTMbije}
Let $X$ be locally compact. 
Let $\mathbf{\Phi^+}  =  \mathbf{\Phi^+} ( C_0^+(X))$.
Consider the map $\Pi : \mathbf{DTM}(X) \longrightarrow \mathbf{\Phi^+} $  where $\Pi(\mu) = \mathcal{R}_{\mu}$ and 
$\mathcal{R}_{\mu}$ is the functional according to Definition \ref{DTM2rho}.
Then the map $\Pi$ has the following properties:
\begin{enumerate}[label=(\Roman*),ref=(\Roman*)]
\item \label{PiConic}
$\Pi$ is conic-linear, i.e.
$\Pi( c \mu  + d \nu) = c \Pi (\mu) +d \Pi(\nu),  \ \ c, d \ge 0.$
\item \label{PiMon} 
$\mu \le \nu $ if and only if $\Pi(\mu) \le \Pi(\nu)$ (i.e., $\mathcal{R}_{\mu} (f)  \le \mathcal{R}_{\nu} (f) $ for all $f \in C_0^+(X)$).
\item \label{TMMbije}
$  \mu \in \mathbf{TM} (X) $ iff $ \rho$ is a quasi-linear functional on $C_0(X)$, and  $\mu \in \mathbf{M}(X)$ iff $\rho$ is a linear functional on $C_0(X)$,
where $ \rho(f) = \Pi(\mu)(f^+) - \Pi(\mu)(f^-) $.
\item \label{normeql}
$\parallel \mathcal{R}_{\mu} \parallel = \mu(X).$
\item
The map $\Pi : \mathbf{DTM} (X) \longrightarrow \mathbf{\Phi^+} $  is a conic-linear order-preserving bijection such that
$\parallel \mathcal{R}_{\mu} \parallel = \mu(X).$
\end{enumerate} 
\end{theorem} 

\begin{proof}
\begin{enumerate}[label=(\Roman*),ref=(\Roman*)]
\item
Let $\lambda = c \mu + d \nu, \, c, d \ge 0$.
Take any $f \in C_0^+(X)$.  For function $R_{1, \lambda, f}$ in Definition \ref{4fns}  
we see that $R_{1, \lambda, f} = c R_{1, \mu, f}  + d R_{1, \nu, f}$.
From formula (\ref{rfformp}) 
$ \mathcal{R}_{\lambda} (f) = c \mathcal{R}_{\mu} (f) + d \mathcal{R}_{\nu} (f)$, and the statement follows.
\item
Let $ \mu \le \nu$. Take any $f \in C_0^+(X)$.  Using Definition \ref{4fns}  we see that $R_{1, \mu, f} \le  R_{1, \nu, f}$. 
Then  by formula (\ref{rfformp}) we have
$\mathcal{R}_{\mu}(f) \le \mathcal{R}_{\nu}(f)$. Thus, $ \mathcal{R}_{\mu} \le \mathcal{R}_{\nu}$.

Now assume that  $ \mathcal{R}_{\mu} \le \mathcal{R}_{\nu}$. From part \ref{MUrhoMU4} of Lemma \ref{MUrhoMU} 
we see that $ \mu(K) \le \nu(K)$ for any compact $K$.
By Remark  \ref{DTMagree} $ \mu \le \nu$.
\item
Suppose $ \rho(f) =  \mathcal{R}_{\mu} (f^+) -  \mathcal{R}_{\mu} (f^-)$ is a quasi-linear functional on $C_0(X)$.
By  part \ref{bijeG} of Corollary \ref{ReprThLC-} and part \ref{QLdastTM} of Theorem \ref{rho2muD} $\mu$ is a finite topological measure. 

Now suppose that $ \mu$ is a finite topological measure.
Let $m_f$ be the measure from Theorem \ref{mf2t} for $f$, where $f \in C_0(X)$.  
Consider functional $ \rho(f) = \int_{\mathbb{R}}  id \,\,  d m_f$ on $C_0(X)$.
For $ \phi \in C(\overline{f(X)}) $ and any open set $W \subseteq \mathbb{R}$  we have
$$ m_{\phi \circ f} (W) = \mu( (\phi \circ f)^{-1} (W) ) = \mu(f^{-1}(\phi^{-1}(W))) = m_f(\phi^{-1}(W)) = (m_f \circ \phi^{-1}) (W),$$
thus, $m_{\phi \circ f}$ and $m_f \circ \phi^{-1}$ are equal as measures on $ \mathbb{R}$, and for $ \phi \circ f \in A(f)$ we obtain
$$ \rho(\phi \circ f) = \int_{\mathbb{R}}  id \, d m_{\phi \circ f} = \int_{\mathbb{R}}  id  \, d(m_f \circ \phi^{-1}) = \int_{\mathbb{R}}  \phi \, dm_f.$$
For  $ \phi \circ f, \psi \circ f \in A(f)$ we have:
$$  \rho(\phi \circ f + \psi \circ f ) = \int_{\mathbb{R}}  (\phi + \psi) \, dm_f = \int_{\mathbb{R}}  \phi  \, dm_f  +  \int_{\mathbb{R}}   \psi \, dm_f 
= \rho(\phi \circ f ) + \rho( \psi \circ f ).$$ 
For any const $c$ we see that $ \rho(cf) = \rho((c \, id) \circ f) = \int_{\mathbb{R}}  c\, id\, dm_f = c \rho(f)$. Thus, $ \rho$ is a quasi-linear functional on $C_0(X)$.
Since $f^{+} = (0 \vee id) \circ f  \in A(f)$, and $ f^{-} \in A(f)$, we know that $\rho(f) = \rho(f^+) - \rho(f^-)$. 
It is clear that $\rho(g) = \mathcal{R}_{\mu} (g) $ for any $ g \in C_0^+$, so  
$\rho(f) = \rho(f^+) - \rho(f^-) =  \mathcal{R}_{\mu} (f^+) -  \mathcal{R}_{\mu} (f^-)  =\Pi(\mu)(f^+) - \Pi(\mu)(f^-) $.

If  $\mu \in \mathbf{M}(X)$ then $ \rho(f) = \int id \,\,  d m_f = \int f \, \, d\mu$ is the usual integral, and the last assertion is a well known fact.
\item
This is part \ref{MUrhoMU6} of Lemma \ref{MUrhoMU}.
\item
Clear from  part \ref{bijeG} of Corollary \ref{ReprThLC-} and parts \ref{PiConic}, \ref{PiMon}, and \ref{normeql}. 
\end{enumerate}
\end{proof}

\begin{theorem} \label{sovpadSem}
Let  $X$  be compact. 
\begin{enumerate}[label=(\roman*),ref=(\roman*)]
\item
$ \mathbf{QI}= \mathbf{\Phi^s}.$
\item
$\mathbf{\Phi^+}  = \mathbf{\Phi^r}$ and $\mathbf{\Phi^-}= \mathbf{\Phi^l} $. 
\end{enumerate} 
Here all functionals are on $ C(X)$. 
\end{theorem}

\begin{proof}
\begin{enumerate}[label=(\roman*),ref=(\roman*)]
\item
By Remark \ref{PhiFam1}  we need to show that $ \mathbf{\Phi^s} \subseteq  \mathbf{QI}$, so let $ \rho \in   \mathbf{\Phi^s}$.  
Then $ \rho \in \mathbf{\Phi^r}$ by Remark \ref{RLD}.  
By Corollary \ref{ReprThLC-} there is a unique deficient topological measure $\mu$ such that 
$\rho (f) = \rho_\mu (f)$ for all $ f \in C(X) $. By part \ref{s1tm} of Theorem \ref{rho2muD}, $\mu$ is a topological measure. 
Then $ \rho$  is a quasi-linear functional (see~\cite[Theorem 4.1]{Aarnes:TheFirstPaper} or~\cite[Theorem 42, Sect. 4]{Butler:QLFLC}).

\item  
If $\rho$ is an r-functional, by Corollary \ref{ReprThLC-} $\rho = \mathcal{R}_{\mu}$ 
for a unique deficient topological measure $\mu$.  
By Theorem \ref{RisA+}  $ \mathcal{R}_{\mu}$ is a  p-conic quasi-linear functional. Thus, $ \mathbf{\Phi^r} \subseteq \mathbf{\Phi^+}  $. 
The other inclusion is given by 
Lemma \ref{+isR}. We can prove that $\mathbf{\Phi^-}= \mathbf{\Phi^l}$ in a similar way, 
using Corollary \ref{ReprThLC-} and Lemma \ref{+isR}. 
\end{enumerate}
\end{proof} 

\begin{remark}
From Theorem \ref{ReprThLC},  part \ref{RTdtmCMP} of Corollary \ref{ReprThLC-}, and Theorem \ref{sovpadSem}, 
it follows that when $X$ is compact, in Theorem \ref{DTMbije} 
we may take $\mathbf{\Phi^+} $ to be $ \mathbf{\Phi^+} (C(X) = \mathbf{\Phi^r}(C(X))$.
\end{remark}

\begin{theorem} \label{Fincls} 
\begin{enumerate}[label=(\Roman*),ref=(\Roman*)]
\item
Let $X$ be locally compact. For functionals on $C_0(X)$ we have:
$$ L \subseteq QI \subseteq \Phi^+ \cap \Phi^- \subseteq  \Phi^r  \cap \Phi^l ,$$
$$ \mathbf{L} \subseteq \mathbf{QI} \subseteq \mathbf{\Phi^+}  \cap \mathbf{\Phi^-} \subseteq  \mathbf{\Phi^r} \cap \mathbf{\Phi^l}   = \mathbf{\Phi^s}. $$
In general, 
$$  L \subsetneqq QI  \subsetneqq \Phi^+. $$
\item  
Let $X$ be compact. Then for functionals on $C(X)$ we have:
$$ \mathbf{L} \subseteq  \mathbf{QI} = \mathbf{\Phi^+}  \cap \mathbf{\Phi^-}=\mathbf{\Phi^r} \cap \mathbf{\Phi^l}  =  \mathbf{\Phi^s}.  $$
In general, 
$$ \mathbf{L} \subsetneqq  \mathbf{QI}. $$
\end{enumerate}
\end{theorem}

\begin{proof}
\begin{enumerate}[label=(\Roman*),ref=(\Roman*)]
\item \label{LCfamil}
The inclusion $QI \subseteq  \Phi^+ \cap \Phi^-$ is given by Remark \ref{QIpnconic}.
The inclusion $\Phi^+ \cap \Phi^- \subseteq  \Phi^r \cap \Phi^l $ follows from Lemma \ref{+isR}.  
By  Remark \ref{RLD} $\Phi^r \subseteq \Phi^c$, so  from Corollary \ref{S=pn}
we see that $ \mathbf{\Phi^r} \cap \mathbf{\Phi^l} \subseteq \mathbf{\Phi^s}$.  By  part \ref{fsINfr} of Remark \ref{RLD} 
we have: $ \Phi^s \subseteq  \Phi^r \cap \Phi^l$. 

The proper inclusion  $L \subsetneqq QI$ follows from the existence of quasi-linear functionals that are not linear, or 
existence of topological measures that are not measures. 
The proper inclusion $QI  \subsetneqq \Phi^+$ follows from  part \ref{TMMbije} of Theorem \ref{DTMbije} 
and the existence of 
deficient topological measures that are not topological measures. See Remark \ref{proinclu}.
For an example of a quasi-linear but not linear functional on 
a locally compact space see~\cite[Example 55, Sect. 5]{Butler:QLFLC}. 
\item  \label{Cfamil}
Use the previous part and Theorem \ref{sovpadSem}.
\end{enumerate}
\end{proof}

\begin{remark} \label{4ways}
Let $X$ be compact.
From Corollary \ref{ReprThLC-} and Theorem \ref{sovpadSem} we see that 
the functionals corresponding to finite deficient topological measures can be described in four ways: 
as p-conic quasi-linear functionals,
as r-functionals, as n-conic quasi-linear functionals, and as l-functionals.  

From Theorem \ref{Fincls} we see that the functionals corresponding to finite topological measures 
can be described in four ways: 
as quasi-linear functionals;   as s-functionals; 
as functionals that are both p-conic quasi-linear  and n-conic quasi-linear; 
and as functionals that are both r- and l-functionals.
\end{remark} 

\begin{remark}
Theorem \ref{Fincls}  answers positively the question posed  in~\cite[Remark 7]{Svistula:DTM}, of whether 
for a compact space $\mathbf{\Phi^s} = \mathbf{\Phi^r} \cap \mathbf{\Phi^l} $.
Note that by part \ref{RLmon} of Lemma \ref{CcondRfl} 
our definition of  r- and l-functionals in the compact case coincide with those in~\cite[Definition 6]{Svistula:DTM}. 
By Definition \ref{sfnl2} and Theorem \ref{Fincls} our definition of an s-functional coincides with the one in~\cite[Definition 6]{Svistula:DTM}, where 
it was first introduced.  
\end{remark}
 
\begin{ack}
This  work was conducted at the Department of Mathematics at the University of California Santa Barbara. 
The author would like to thank the department for its hospitality and supportive environment.
\end{ack}


\small{

\bibliographystyle{abbrv}


\begin{thebibliography}{1}

\bibitem{Aarnes:TheFirstPaper}
 J. ~Aarnes.
 \newblock Quasi-states and quasi-measures. 
 \newblock {\em Adv. Math.}, 86 (1): 41--67, 1991.   
      
\bibitem{Butler:TechniqLC}
  S. ~Butler. 
  \newblock Ways of obtaining topological measures on locally compact spaces.
  \newblock {\em Bulletin of Irkutsk State University, Series "Mathematics"}, 25: 33--45, 2018.
  
\bibitem{Butler:TMLCconstr}
  S. ~Butler. 
  \newblock Solid-set functions and topological measures on locally compact spaces. arXiv: 1902.01957  
      
\bibitem{Butler:DTMLC}
  S. ~Butler. 
  \newblock Deficient topological measures on locally compact spaces. arXiv: 1902.02458   
  
\bibitem{Butler:QLFLC}
  S. ~Butler.
  \newblock Quasi-linear functionals on locally compact spaces.  arXiv: 1902.03358 
  
\bibitem{Dugundji}
J. ~Dugundji.
\newblock {\em Topology}.
\newblock  Allyn and Bacon, Inc., Boston, 1966.

\bibitem{Engelking}
R. ~Engelking. 
\newblock {\em General topology}.
\newblock PWN, Warsaw, 1989.

\bibitem{Grubb:Signed} 
  D. ~Grubb. 
 \newblock Signed quasi-measures.
 \newblock {\em Trans. Amer. Math. Soc.}, 349(3): 1081--1089, 1997.
 
\bibitem{Grubb:Lectures} 
 D. ~Grubb. 
 \newblock Lectures on quasi-measures and quasi-linear functionals on compact spaces. Unpublished, 1998.
  
\bibitem{OrjanAlf:CostrPropQlf}
  $\O$. ~Johansen,  A. ~Rustad.
 \newblock Construction and Properties of quasi-linear functionals.
 \newblock{\em Trans. Amer. Math. Soc.} 358(6): 2735--2758, 2006.  
  
\bibitem{EntovPolterovich}
M. ~Entov, L. Polterovich. 
 \newblock Quasi-states and symplectic intersections. 
 \newblock{ \em Comment. Math. Helv.}, 81: 75--99, 2006.

\bibitem{PoltRosenBook}
L. ~Polterovich, D. ~Rosen.
\newblock {\em Function theory on symplectic manifolds}.
\newblock CRM Monograph series, vol. 34.  American Mathematical Society, Providence,
Rhode Island, 2014.
  
\bibitem{Alf:ReprTh}
  A. ~Rustad. 
  \newblock Unbounded quasi-integrals.
  \newblock {\em Proc. Amer. Math. Soc.}, 129(1): 165--172, 2000.    

\bibitem{Svistula:Signed}
M. ~Svistula. 
\newblock A Signed quasi-measure decomposition.
\newblock {\em Vestnik Samara Gos. Univ. Estestvennonauchn.}, 62(3):192--207, 2008. (in Russian)

\bibitem{Svistula:DTM}
M. ~Svistula.
\newblock Deficient topological measures and functionals generated by them.
\newblock {\em Sbornik: Mathematics}, 204(5): 726--76, 2013.
  
\end{thebibliography}

}

\end{document}